\documentclass[10pt,twocolumn,twoside]{IEEEtran}              
\IEEEoverridecommandlockouts         

\usepackage{amsmath,amsfonts,amssymb,color}
\usepackage{epsfig}
\usepackage{psfrag}
\usepackage{algorithm}
\usepackage[noend]{algpseudocode}
\usepackage{array}
\usepackage{epstopdf}
\usepackage{cite}
\usepackage{footmisc}
\usepackage{mathtools}
\usepackage{graphicx}
\usepackage{xcolor}
\usepackage{caption}
\usepackage{subfig}

\newtheorem{problem}{Problem}
\newtheorem{theorem}{Theorem}

\newtheorem{lemma}{Lemma}
\newtheorem{remark}{Remark}
\newtheorem{definition}{Definition}

\newtheorem{example}{Example}
\newtheorem{assumption}{Assumption}
\renewcommand{\theenumi}{(\alph{enumi}}
\algnewcommand{\LineComment}[1]{\State \(\triangleright\) #1}

\title{\LARGE \bf
Resilient Sensor Placement for Kalman Filtering in Networked Systems: Complexity and Algorithms
}
\begin{document}

\author{Lintao Ye, Sandip Roy and Shreyas Sundaram
\thanks{This research was supported by NSF grants CMMI-1635014 and CMMI-1635184.  Lintao Ye and Shreyas Sundaram are with the School of Electrical
and Computer Engineering at Purdue University, West Lafayette, IN 47907 USA. Email: \{ye159,sundara2\}@purdue.edu.  Sandip Roy is with the School of Electrical Engineering and Computer Science at Washington State University, Pullman, WA 99164 USA.  Email: sroy@eecs.wsu.edu.
}
}
\maketitle
\begin{abstract}
Given a linear dynamical system affected by noise, we study the problem of optimally placing sensors (at design-time) subject to a sensor placement budget constraint in order to minimize the trace of the steady-state error covariance of the corresponding Kalman filter. While this problem is NP-hard in general, we consider the underlying graph associated with the system dynamics matrix, and focus on the case when there is a single input at one of the nodes in the graph. We provide an optimal strategy (computed in polynomial-time) to place the sensors over the network. Next, we consider the problem of attacking (i.e., removing) the placed sensors under a sensor attack budget constraint in order to maximize the trace of the steady-state error covariance of the resulting Kalman filter. Using the insights obtained for the sensor placement problem, we provide an optimal strategy (computed in polynomial-time) to attack the placed sensors. Finally, we consider the scenario where a system designer places the sensors under a sensor placement budget constraint, and an adversary then attacks the placed sensors subject to a sensor attack budget constraint. The resilient sensor placement problem is to find a sensor placement strategy to minimize the trace of the steady-state error covariance of the Kalman filter corresponding to the sensors that survive the attack. We show that this problem is NP-hard, and provide a pseudo-polynomial-time algorithm to solve it.
\end{abstract}

\section{Introduction}
In large-scale control system design, one of the key problems is to place sensors or actuators on the system in order to achieve certain performance criteria (e.g., \cite{olshevsky2014minimal}, \cite{pequito2016minimum}). In cases involving linear systems with process or measurement noise, researchers have studied how to place sensors (at design-time) in order to minimize certain metrics of the error covariance of the corresponding Kalman filter (e.g., \cite{chmielewski2002theory,mo2011sensor,yang2015deterministic,tzoumas2016sensor,Zhang2017sensor}). The problem has been shown to be NP-hard and inapproximable within any constant factor in general  \cite{ye2018complexity}. This motivates us to consider special classes of this problem in this paper and seek polynomial-time algorithms for the optimal sensor placement problem. Specifically, we consider a discrete-time linear dynamical system whose states represent nodes in a directed network, and interact according to the topology of the network. The nodes of the network are possibly affected by stochastic inputs. Such networked systems with stochastic inputs have received much attention from researchers recently (e.g., \cite{xiao2007distributed,lin2014algorithms,clark2014supermodular,fitch2016joint}). These models encompass diffusion networks (e.g., \cite{koorehdavoudi2017interactions}) that arise in many different areas, including information and influence diffusion over social networks \cite{kempe2003maximizing}, spreading of diseases in populations \cite{goltsev2012localization}, and diffusion of chemicals in certain environments \cite{roy2015situational}. In such applications, estimating the states of the entire network is an important objective.

In this paper, we focus on the case when there is a single node of the network that has a stochastic input. Specifically, we consider the scenario where a system designer can spend a limited budget on placing sensors (at design-time) over the network in order to minimize the trace of the steady-state error covariance of the Kalman filter corresponding to the placed sensors. A sensor placed at a certain node gives measurements of the state corresponding to the node. In addition, placing a sensor at a node incurs a placement cost (which could vary across the nodes). We refer to this problem as the {\it Graph-based Kalman Filtering Sensor Placement (GKFSP)} problem. 

Additionally, the systems that we are interested in monitoring may be targeted by adversaries, where an adversary can attack a subset of placed sensors. Different types of attacks have been studied previously, including Denial of Service (DoS) attacks (e.g., \cite{wood2002denial}, \cite{zhang2015optimal}) and false data injection attacks (e.g., \cite{mo2010false}, \cite{liu2011false}). Here, we consider adversaries that perform DoS attacks on sensors by simply removing them (or equivalently, dropping all the measurement data). The goal of the adversary is to remove a subset of placed sensors under a budget constraint in order to maximize the trace of the steady-state error covariance of the Kalman filter corresponding to the surviving sensors. We assume that attacking a sensor placed at a node incurs an attack cost (which could also vary across the nodes). In contrast with existing work in the literature, we analyze the problem using the graph structure of the systems.  We refer to this problem as the {\it Graph-based Kalman Filtering Sensor Attack (GKFSA)} problem. 

Finally, combining the two problems that we considered above, we formulate and study a resilient sensor placement problem for the networked system. We assume that the system designer is aware of the potential attack from an adversary who chooses to optimally attack the sensors (subject to an attack budget constraint) deployed by the system designer. The system designer's goal is to place sensors (under a placement budget constraint) among a subset of nodes in order to minimize the trace of the steady-state error covariance of the Kalman filter corresponding to the surviving sensors after the attack. We refer to this problem as the  {\it Resilient Graph-based  Kalman Filtering Sensor Placement (RGKFSP)} problem. 
\vspace{-0.2cm}
\subsection*{Contributions}
First, we provide an optimal sensor placement strategy for the GKFSP problem using the graph structure of the system. Second, leveraging the insights for the GKFSP problem, we give an optimal sensor attack strategy for the GKFSA problem. Third, we show that the RGKFSP problem is NP-hard; we then provide an algorithm based on dynamic programming that can return an optimal solution to general instances of the RGKFSP problem in pseudo-polynomial time. Although the results are derived under the assumption that the sensors give perfect measurements, we show that how to apply these results to analyze the case with sensor measurement noise and provide numerical examples. A preliminary version of the results in this paper was presented in \cite{ye2018optimal}, where only the GKFSP problem was studied for a more restrictive class of system dynamics matrices.
\vspace{-0.2cm}
\subsection*{Related work}
The (design-time) sensor placement problem has been widely studied in the literature. For example, in \cite{shamaiah2010greedy,tzoumas2016sensor}, the authors considered the Kalman filtering sensor placement problem over a finite number of time steps. Here, we study the problem of optimizing steady-state error covariances of the corresponding Kalman filter.  In \cite{Zhang2017sensor,ye2018complexity}, the authors considered the same sensor placement problem as the one considered here, but for general system dynamics.  In such cases, they showed that finding the optimal placement for the general problem is NP-hard.  Thus, in this paper, we impose additional structure on the problem (by considering the graph representation of the dynamics) in order to seek optimal solutions. In \cite{liu2016sensor,joshi2008sensor}, the authors studied the sensor placement problem for estimating a static variable (parameter) that does not change over time. Here, we study the problem of placing sensors to estimate the states of a linear dynamical system affected by stochastic inputs.  In contrast to the sensor placement problem where the set of placed sensors cannot change over time, the sensor scheduling problem for Kalman filtering has also received much attention (e.g., \cite{gupta2006stochastic,vitus2012efficient,huber2012optimal,jawaid2015submodularity}), where different sets of sensors can be chosen at different time steps.

In networked system settings, the authors in \cite{roy2018graph} considered the sensor placement problem for continuous-time diffusion dynamics, and applied the Wiener filter to estimate the system states using sensor measurements. Here, we consider discrete-time networked system dynamics and apply the Kalman filter to estimate the system states. The authors in \cite{lin2014algorithms,clark2014supermodular,fitch2016joint} studied the leader selection problem in consensus networks with stochastic inputs. The problem is to select a subset of nodes whose states are fixed over time in order to minimize the $H_2$ norm of the system states at steady state. In contrast, we consider the problem of placing sensors among the nodes of systems with more general dynamics in order to minimize the trace of the steady-state error covariance of the Kalman filter.

Although both of the sensor placement and the sensor attack problems have received much attention from researchers, the resilient sensor placement is less explored.  
The authors in \cite{tzoumas2017resilient} considered the problem of resilient maximization of monotone submodular set functions under a cardinality constraint on the sets. They proposed a polynomial-time approximation algorithm for the problem with performance bounds that depend on the curvature of the objective function. In \cite{laszka2015resilient}, the authors considered a resilient observation selection problem. The problem is to resiliently select observations of a scalar Gaussian process given that some of the selected observations could be removed by an adversary. The authors showed that this problem is NP-hard and proposed a greedy algorithm with a provable performance guarantee. Here, we consider the resilient sensor placement problem for Kalman filtering of (vector) linear dynamical systems subject to general knapsack constraints. While we show this problem is NP-hard, we give an algorithm based on dynamic programming to solve the problem optimally in pseudo-polynomial time \cite{garey1979computers}.

\subsection*{Notation and terminology}
The sets of integers and real numbers are denoted as $\mathbb{Z}$ and $\mathbb{R}$, respectively. For any $x\in\mathbb{R}$, let $\lfloor x\rfloor$ denote the greatest integer that is less than or equal to $x$. For a matrix $P\in\mathbb{R}^{n\times n}$, let $P^T$ denote its transpose, $P_{ij}$ (or $(P)_{ij}$) denote the element in the $i$th row and $j$th column of $P$, and $P_i$ denote the $i$th row of $P$. Let $\mathbf{0}_{m\times n}$ denote a zero matrix; the subscript is dropped if the dimension of the matrix is clear from the context. The identity matrix of dimension $n$ is denoted as $I_n$. A positive semi-definite matrix $P$ is denoted by $P\succeq\mathbf{0}$ and $P\succeq Q$ if $P-Q\succeq\mathbf{0}$. The set of $n$ by $n$ positive definite (resp., positive semi-definite) matrices is denoted by $\mathbb{S}_{++}^n$ (resp., $\mathbb{S}_+^n$).  For a vector $x$, denote its $i$th element as $x_i$, and let $\textrm{supp}(x)$ be its support, where $\textrm{supp}(x)=\{i:x_i\neq 0\}$. Define $\mathbf{e}_i$ to be a column vector where the $i$th element is $1$ and all the other elements are zero; the dimension of the vector can be inferred from the context. We use $\mathbb{E}[x]$ to denote the expectation of a random variable (vector) $x$. For a set $\mathcal{A}$, let $|\mathcal{A}|$ be its cardinality. Given two functions $\varphi_1:\mathbb{R}_{\ge0}\to\mathbb{R}$ and $\varphi_2:\mathbb{R}_{\ge0}\to\mathbb{R}$, $\varphi_1(n)$ is $O(\varphi_2(n))$ if there exist positive constants $c$ and $N$ such that $|\varphi_1(n)|\le c|\varphi_2(n)|$ for all $n\ge N$.

\section{Problem Formulation} \label{sec:problem formulation}
We begin with the following definitions from graph theory. Further details can be found in, for example, \cite{bondy1976graph} and \cite{horn1985matrix}.
\begin{definition}
\label{def:the associated graph of a matrix}
For any given matrix $A\in\mathbb{R}^{n\times n}$, the directed graph of $A$, denoted as $\mathcal{G}(A)$, is defined as the directed graph on $n$ vertices (or nodes) $x_1,x_2,\dots,x_n$ such that for all $i,j\in\{1,2,\dots,n\}$, there is a directed edge in $\mathcal{G}(A)$ from $x_j$ to $x_i$, denoted as $(x_j,x_i)$, if and only if $A_{ij}\neq 0$. Denoting the set of vertices and the set of edges of $\mathcal{G}(A)$ as $\mathcal{X}(A)\triangleq\{x_1,x_2,\dots,x_n\}$ and $\mathcal{E}(A)$, respectively, the graph $\mathcal{G}(A)$ is also denoted as $\mathcal{G}(A)=\{\mathcal{X}(A),\mathcal{E}(A)\}$.
\end{definition}

\begin{definition}
\label{def:directed path}
Consider a directed graph $\mathcal{G}=\{\mathcal{X},\mathcal{E}\}$, where $\mathcal{X}\triangleq\{x_1,x_2,\dots,x_n\}$. A directed path from $x_{i_0}$ to $x_{i_t}$ is a sequence of directed edges $(x_{i_0},x_{i_1}),(x_{i_1},x_{i_2}),\dots,(x_{i_{t-1}},x_{i_t})$ in $\mathcal{G}$. The ordered list of vertices in the directed path is $x_{i_0},x_{i_1},\dots,x_{i_t}$. The length of a directed path is the number of directed edges in the directed path. A cycle is a directed path that begins and ends at the same vertex which occurs exactly twice in the ordered list of vertices in the directed path, and no other vertices occur more than once in the list. A cycle of length $1$ is a self-loop at the corresponding vertex.
\end{definition}

\begin{definition}
Consider a directed graph $\mathcal{G}=\{\mathcal{X},\mathcal{E}\}$. For any pair of distinct vertices $x_i,x_j\in\mathcal{X}$ such that there exists a directed path from $x_i$ to $x_j$, the distance from $x_i$ to $x_j$, denoted as $l_{ij}$, is defined as the shortest length over all such paths. Define $l_{mm}=0$ for all $x_m\in\mathcal{X}$.
\end{definition}

\begin{definition}
\label{def:strongly connected}
A directed graph $\mathcal{G}=\{\mathcal{X},\mathcal{E}\}$ is strongly connected if for all pairs of distinct vertices $x_i,x_j\in\mathcal{X}$, there is a directed path from $x_j$ to $x_i$ in $\mathcal{G}$. 
\end{definition}

We start with a general system model. Consider a matrix $A\in\mathbb{R}^{n\times n}$ with the associated graph $\mathcal{G}(A)=\{\mathcal{X}(A),\mathcal{E}(A)\}$ (given in Definition $\ref{def:the associated graph of a matrix}$). Suppose that $\mathcal{I}\triangleq\{x_{i_0},x_{i_1},\dots,x_{i_{n_1-1}}\}\subseteq\mathcal{X}(A)$ is the set of nodes that have stochastic inputs, where $n_1\in\mathbb{Z}_{\ge1}$. We then consider the following discrete-time linear system:
\begin{equation}
x[k+1] = Ax[k] + Bw[k],
\label{eqn:system dynamics}
\end{equation}
where $x[k]\in\mathbb{R}^n$ is the system state at time step $k$, and $B\triangleq\begin{bmatrix}\mathbf{e}_{i_0} & \cdots & \mathbf{e}_{i_{n_1-1}}\end{bmatrix}\in\mathbb{R}^{n\times n_1}$ is the input matrix. The stochastic input $w[k]\in\mathbb{R}^{n_1}$ is a zero-mean white noise process with $\mathbb{E}[w[k](w[k])^T]=W\in\mathbb{S}_+^{n_1}$. The initial state $x[0]$ is a random vector with mean $\bar{x}_0\in\mathbb{R}^n$ and covariance $\Pi_0\in\mathbb{S}_+^n$, and is assumed to be independent of $w[k]$ for all $k\in\mathbb{Z}_{\ge0}$. Each state of the system, denoted as $x_i[k]$, is associated with  node $x_i$ in $\mathcal{G}(A)$. As we mentioned in the introduction, \cite{ye2018complexity} showed that the Kalman filtering sensor placement problem cannot be approximated within any constant factor in polynomial time (if P$\neq$NP) for general system dynamics matrices {\it even} when the measurement noise is zero. Moreover, under the networked system setting, \cite{ye2018optimal} showed that if there are multiple input nodes in the graph, the Kalman filtering sensor placement problem becomes NP-hard {\it even} when the graph only contains a set of disjoint paths of length three and each path has a single input node. Hence, in order to bypass these inherent complexity issues, we focus on networked systems with a single input node $x_{i_0}\in\mathcal{X}(A)$ (i.e., $B=\mathbf{e}_{i_0}$ and $\mathbb{E}[(w[k])^2]=\sigma^2_w\in\mathbb{R}_{\ge0}$), and  seek efficient algorithms to optimally solve the corresponding sensor placement, sensor attack, and resilient sensor placement problems.  We assume throughout this paper that the pair $(A,B\sigma_w)$ is stabilizable. The generality of this assumption will be justified later.

\subsection{The Sensor Placement Problem}
First, suppose that there is a system designer who can choose a subset of the vertices of the graph $\mathcal{G}(A)$ at which to place sensors under a budget constraint. Specifically, a sensor placed at node $x_i\in\mathcal{X}(A)$ has a placement cost $h_{i}\in\mathbb{Z}_{\ge0}$; define the sensor placement cost vector as $h\triangleq \left[\begin{matrix}h_1 & \cdots & h_n\end{matrix}\right]^T$. The designer has a sensor placement budget $H\in\mathbb{Z}_{\ge0}$ that can be spent on placing sensors at the nodes of $\mathcal{G}(A)$. A sensor that is placed at node $x_i\in\mathcal{X}(A)$ gives a measurement
\begin{equation}
\label{eqn:single sensor measurement}
y_i[k]=C_ix[k]+v_i[k],
\end{equation}
where $C_i=\mathbf{e}_i^T$ and $v_i[k]\in\mathbb{R}$ is a zero-mean white noise process. We further define $y[k]\triangleq\big[y_1[k]\ \cdots \ y_n[k]\big]^T$, $C\triangleq\big[C_1^T\ \cdots\ C_n^T\big]^T$ and $v[k]\triangleq\big[v_1[k]\ \cdots \ v_n[k]\big]^T$. Thus, the output provided by all sensors together is given by
\begin{equation}
\label{eqn:all sensors measurements}
y[k]=Cx[k]+v[k],
\end{equation}
where $C=I_{n}$. We denote $\mathbb{E}[v[k](v[k])^T]=V\in\mathbb{S}_{+}^n$ and consider $\mathbb{E}[v[k](w[j])^T]=\mathbf{0}$, $\forall k, j\in\mathbb{Z}_{\ge0}$. The initial state $x[0]$ is also assumed to be independent of $v[k]$ for all $k\in\mathbb{Z}_{\ge0}$.

After the sensors are placed, the Kalman filter is then applied to provide an estimate of the states using the measurements from the installed sensors. We define a vector $\mu\in\{0,1\}^n$ as the indicator vector indicating the vertices where sensors are placed. Specifically, $\mu_i=1$ if and only if a sensor is placed at node $x_i\in\mathcal{X}(A)$. Denote $C(\mu)$ as the measurement matrix of the installed sensors indicated by $\mu$, i.e., $C(\mu)\triangleq\left[\begin{matrix} C_{i_1}^T &  \cdots & C_{i_p}^T \end{matrix}\right]^T$, where $\textrm{supp}(\mu)=\{i_1,\dots,i_p\}\subseteq\{1,\dots,n\}$. Similarly, denote $V(\mu)$ as the measurement noise covariance matrix of the installed sensors, i.e., $V(\mu)=\mathbb{E}[\tilde{v}[k](\tilde{v}[k])^T]$, where $\tilde{v}[k]\triangleq\big[(v[k])_{i_1}\ \cdots \ (v[k])_{i_p}\big]^T$. The {\it a priori} and the {\it a posteriori} error covariance matrices of the Kalman filter at time step $k$, when the sensors indicated by $\mu$ are placed, are denoted as $\Sigma_{k/k-1}(\mu)$ and $\Sigma_{k/k}(\mu)$, respectively. The initial {\it a priori} error covariance is set as $\Sigma_{0/-1}(\mu)=\Pi_0$.
The limit $\Sigma(\mu)\triangleq\mathop{\lim}_{k\to\infty}\Sigma_{k+1/k}$ (also known as the steady-state {\it a priori} error covariance), if it exists, satisfies the \textit{discrete algebraic Riccati equation (DARE)} \cite{anderson1979optimal}:
\begin{multline}
\label{eqn:DARE}
\Sigma(\mu)=A\Sigma(\mu)A^T+\sigma_w^2BB^T - \\
A\Sigma(\mu)C(\mu)^T\big(C(\mu)\Sigma(\mu)C(\mu)^T+V(\mu)\big)^{-1}C(\mu)\Sigma(\mu)A^T,
\end{multline}
where $\sigma_w^2\in\mathbb{R}_{\ge0}$ and $B=\mathbf{e}_{i_0}$. The limit $\Sigma^*(\mu)\triangleq\mathop{\lim}_{k\to\infty}\Sigma_{k/k}(\mu)$ (also known as the steady-state {\it a posteriori} error covariance), if it exists, satisfies the following equations \cite{catlin1989estimation}:
\begin{multline}
\label{eqn:postDARE}
\Sigma^*(\mu)=\Sigma(\mu) -\\ \Sigma(\mu)C(\mu)^T\big(C(\mu)\Sigma(\mu)C(\mu)^T+V(\mu)\big)^{-1}C(\mu)\Sigma(\mu),
\end{multline}
and 
\begin{equation}
\label{eqn:coupled DARE}
\Sigma(\mu)=A\Sigma^*(\mu)A^T+\sigma_w^2BB^T.
\end{equation}
The inverses in Eq. \eqref{eqn:DARE} and \eqref{eqn:postDARE} are interpreted as the Moore-Penrose pseudo-inverses (which we denote using the notation ``$\dagger$") if the arguments are not invertible \cite{anderson1979optimal}. We will use the following result from \cite{anderson1979optimal}.
\begin{lemma}
For a given indicator vector $\mu$, $\Sigma_{k/k-1}(\mu)$ (resp., $\Sigma_{k/k}(\mu)$) will converge, as $k\to\infty$, to a finite limit $\Sigma(\mu)$ (resp., $\Sigma^*(\mu)$), regardless of the initial covariance $\Sigma_{0/-1}(\mu)$, if and only if the pair $(A,C(\mu))$ is detectable and the pair $(A,B\sigma_w)$ is stabilizable. Furthermore, if the limit $\Sigma(\mu)$ (resp., $\Sigma^*(\mu)$) exists, it is also the only positive semi-definite solution to Eq. \eqref{eqn:DARE} (resp., Eq. \eqref{eqn:postDARE}).
\label{lemma:Anderson optimal filtering graph}\label{lemma:Anderson optimal filtering}
\end{lemma}

When the pair $(A,C(\mu))$ is not detectable, we define the limits $\Sigma(\mu)=+\infty$ and $\Sigma^*(\mu)=+\infty$. The priori and posteriori Graph-based Kalman Filtering Sensor Placement (GKFSP) problems are defined as follows.
\begin{problem}
\label{pro:GKFSP problem}
(Priori and Posteriori GKFSP) Consider a system dynamics matrix $A\in\mathbb{R}^{n\times n}$ with the associated graph $\mathcal{G}(A)=\{\mathcal{X}(A),\mathcal{E}(A)\}$, a single vertex $x_{i_0}\in\mathcal{X}(A)$ that has a stochastic input with variance $\sigma_w^2\in\mathbb{R}_{\ge0}$,  the measurement matrix $C=I_n$ (containing all of the individual sensor measurement matrices), a sensor noise covariance matrix $V\in\mathbb{S}_+^n$, a sensor placement cost vector $h\in\mathbb{Z}^n_{\ge0}$ and a sensor placement budget $H\in\mathbb{Z}_{\ge0}$. The priori Graph-based Kalman Filtering Sensor Placement (GKFSP) problem is to find the sensor placement $\mu$, i.e., the indicator vector $\mu$ of the vertices where sensors are placed, that solves
\begin{equation*}
\begin{split}
&\mathop{\min}_{\mu\in\{0,1\}^n} \text{trace}(\Sigma(\mu))\\
&s.t.\ h^T \mu\le H,
\end{split}
\end{equation*}
where $\Sigma(\mu)$ is given by Eq. \eqref{eqn:DARE} if the pair $(A,C(\mu))$ is detectable, and $\Sigma(\mu)=+\infty$ otherwise.  The posteriori GKFSP Problem is to find the sensor placement $\mu$ that solves
\begin{equation*}
\begin{split}
&\mathop{\min}_{\mu\in\{0,1\}^n} \text{trace}(\Sigma^*(\mu))\\
&s.t.\ h^T \mu\le H,
\end{split}
\end{equation*}
where $\Sigma^*(\mu)$ is given by Eq. \eqref{eqn:postDARE} if the pair $(A,C(\mu))$ is detectable, and $\Sigma^*(\mu)=+\infty$ otherwise.
\end{problem}

\subsection{The Sensor Attack Problem}
Suppose that the sensors indicated by the sensor placement $\mu\in\{0,1\}^n$ are placed and installed by the system designer, and there is an adversary who aims to attack (i.e., remove) a subset of the installed sensors. To attack a sensor placed at node $x_i\in\mathcal{X}(A)$, the adversary needs to pay a cost $f_i\in\mathbb{Z}_{\ge0}$. Define the sensor attack cost vector as $f\triangleq\begin{bmatrix}f_1 & \cdots & f_n\end{bmatrix}^T$. The adversary has a total sensor attack budget $F\in\mathbb{Z}_{\ge0}$ for attacking the installed sensors. We define a vector $\nu\in\{0,1\}^n$ as the indicator vector indicating the subset of sensors that are attacked, where $\nu_i=1$ if and only if the sensor at $x_i\in\mathcal{X}(A)$ is attacked. Note that $\text{supp}(\nu)\subseteq\text{supp}(\mu)$ is always assumed implicitly in the sequel. Denote the matrix $C(\mu\setminus\nu)$ as the measurement matrix of the surviving sensors corresponding to $\mu$ and $\nu$, i.e., $C(\mu\setminus\nu)\triangleq\begin{bmatrix}C^T_{j_1} & \cdots & C^T_{j_q}\end{bmatrix}^T$, where $\{j_1,\dots,j_q\}=\text{supp}(\mu)\setminus\text{supp}(\nu)$. Denote $\text{supp}(\mu)\setminus\text{supp}(\nu)\triangleq\text{supp}(\mu\setminus\nu)$. Similarly, define $V(\mu\setminus\nu)$ as the measurement noise covariance of the surviving sensors. The Kalman filter is then applied based on the measurements of the surviving sensors. The resulting {\it a priori} and {\it a posteriori} error covariances of the Kalman filter at time step $k$ are denoted as $\Sigma_{k/k-1}(\mu\setminus\nu)$ and $\Sigma_{k/k}(\mu\setminus\nu)$, respectively, whose limits as $k\to\infty$ are denoted as $\Sigma(\mu\setminus\nu)$ and $\Sigma^*(\mu\setminus\nu)$, respectively. 

The priori and posteriori Graph-based Kalman Filtering Sensor Attack (GKFSA) problems are then defined as follows.
\begin{problem}
\label{pro:GKFSA problem}
(Priori and Posteriori GKFSA) Consider a system dynamics matrix $A\in\mathbb{R}^{n\times n}$ with the associated graph $\mathcal{G}(A)=\{\mathcal{X}(A),\mathcal{E}(A)\}$, a single vertex $x_{i_0}\in\mathcal{X}(A)$ that has a stochastic input with variance $\sigma_w^2\in\mathbb{R}_{\ge0}$, the measurement matrix $C=I_n$ (containing all of the individual sensor measurement matrices), a sensor noise covariance matrix $V\in\mathbb{S}_+^n$, a sensor attack cost vector $f\in\mathbb{Z}^n_{\ge0}$, a sensor attack budget $F\in\mathbb{Z}_{\ge0}$, and a sensor placement vector $\mu\in\{0,1\}^n$. The priori Graph-based Kalman Filtering Sensor Attack (GKFSA) problem is to find the sensor attack $\nu$, i.e., the indicator vector $\nu$ of the vertices where the installed sensors (indicated by $\mu$) are attacked, that solves
\begin{equation*}
\begin{split}
&\mathop{\max}_{\nu\in\{0,1\}^n} \text{trace}(\Sigma(\mu\setminus\nu))\\
&s.t.\ f^T \nu\le F,
\end{split}
\end{equation*}
where $\Sigma(\mu\setminus\nu)$ is given by Eq. \eqref{eqn:DARE} if the pair $(A,C(\mu\setminus\nu))$ is detectable, and $\Sigma(\mu\setminus\nu)=+\infty$ otherwise. The posteriori GKFSA problem is to find the sensor attack $\nu$ that solves
\begin{equation*}
\begin{split}
&\mathop{\max}_{\nu\in\{0,1\}^n} \text{trace}(\Sigma^*(\mu\setminus\nu))\\
&s.t.\ f^T \nu\le F,
\end{split}
\end{equation*}
where $\Sigma^*(\mu\setminus\nu)$ is given by Eq. \eqref{eqn:postDARE} if the pair $(A,C(\mu\setminus\nu))$ is detectable, and $\Sigma^*(\mu\setminus\nu)=+\infty$ otherwise.
\end{problem}

\subsection{The Resilient Sensor Placement Problem}
We next consider the scenario where the system designer is aware of the potential attack from a strategic adversary (who can perform optimal sensor attacks under budget constraints), and aims to choose a resilient sensor placement.We first define feasible sensor placements for the system designer as follows.
\begin{definition}
\label{def:feasible sensor placement}
A sensor placement $\mu\in\{0,1\}^n$ is said to be feasible if $h^T\mu\le H$ (i.e., the sensor placement budget constraint is satisfied), and for all $\nu\in\{0,1\}^n$ such that  $f^T\nu\le F$, $\text{supp}(\mu\setminus\nu)\neq\emptyset$ (i.e., for all sensor attacks that satisfy the sensor attack budget constraint, at least one sensor indicated by $\mu$ is left over by the adversary).
\end{definition} 

\begin{remark}
Note that if a sensor placement $\mu$ is not feasible, there is an attack (satisfying the attacker's budget constraint) such that that the pair $(A,C(\mu\setminus\nu))$ is not detectable if the system dynamics matrix $A$ is not stable.
\end{remark}

The priori and posteriori Resilient Graph-based Kalman Filtering Sensor Placement (RGKFSP) problems are then given by the following.
\begin{problem}
\label{pro:RGKFSP problem}
(Priori and Posteriori RGKFSP) Consider a system dynamics matrix $A\in\mathbb{R}^{n\times n}$ with the associated graph $\mathcal{G}(A)=\{\mathcal{X}(A),\mathcal{E}(A)\}$, a single vertex $x_{i_0}\in\mathcal{X}(A)$ that has a stochastic input with variance $\sigma_w^2\in\mathbb{R}_{\ge0}$, the measurement matrix $C=I_n$ (containing all of the individual sensor measurement matrices), a sensor noise covariance matrix $V\in\mathbb{S}_+^n$, a sensor placement cost vector $h\in\mathbb{Z}^n_{\ge0}$, a sensor placement budget $H\in\mathbb{Z}_{\ge0}$, a sensor attack cost vector $f\in\mathbb{Z}^n_{\ge0}$, and a sensor attack budget $F\in\mathbb{Z}_{\ge0}$. The priori Resilient Graph-based Kalman Filtering Sensor Placement (RGKFSP) problem is to find the sensor placement $\mu$ that solves
\begin{equation*}
\begin{split}
&\mathop{\min}_{\mu\in\{0,1\}^n}\mathop{\max}_{\nu\in\{0,1\}^n} \text{trace}(\Sigma(\mu\setminus\nu)) \\
&s.t.\ h^T \mu\le H,\ \text{and}\ f^T\nu\le F,
\end{split}
\end{equation*}
where $\Sigma(\mu\setminus\nu)$ is given by Eq. \eqref{eqn:DARE} if the pair $(A,C(\mu\setminus\nu))$ is detectable, and $\Sigma(\mu\setminus\nu)=+\infty$ otherwise. The posteriori RGFKSP problem is to find the sensor placement $\mu$ that solves
\begin{equation*}
\begin{split}
&\mathop{\min}_{\mu\in\{0,1\}^n}\mathop{\max}_{\nu\in\{0,1\}^n} \text{trace}(\Sigma^*(\mu\setminus\nu)) \\
&s.t.\ h^T \mu\le H,\ \text{and}\ f^T\nu\le F,
\end{split}
\end{equation*}
where $\Sigma^*(\mu\setminus\nu)$ is given by Eq. \eqref{eqn:postDARE} if the pair $(A,C(\mu\setminus\nu))$ is detectable, and $\Sigma^*(\mu\setminus\nu)=+\infty$ otherwise.
\end{problem}

\section{Solving the GKFSP and GKFSA problems}\label{sec:alg for KFSP and KFSA}
In this section, we provide algorithms to optimally solve the GKFSP and GKFSA problems, respectively, when the sensor noise covariance is $V=\mathbf{0}_{n\times n}$. We will make the following assumptions on the instances of the GKFSP and GKFSA problems in the sequel.
\begin{assumption}
\label{assumption:stabilizable and detectable}
The pair $(A,B\sigma_w)$ is assumed to be stabilizable. The pair $(A,C(\mu))$ is assumed to be detectable for all sensor placements $\mu\in\{0,1\}^n$ with $\text{supp}(\mu)\neq\emptyset$.
\end{assumption}
\begin{assumption}
\label{assumption:distance}
The graph $\mathcal{G}(A)=\{\mathcal{X}(A),\mathcal{E}(A)\}$ (associated with the system dynamics matrix $A\in\mathbb{R}^{n\times n}$) is assumed to satisfy the property that for all $x_j\in\mathcal{X}(A)$ and $x_j\neq x_{i_0}$,  there exists a directed path from $x_{i_0}$ to $x_j$. The system dynamics matrix $A$ is assumed to satisfy $(A^m)_{ji_0}\neq0$ if $l_{i_0j}=m$, where $l_{i_0j}$ is the distance from $x_{i_0}$ to $x_j$.
\end{assumption}
\begin{remark}
\label{remark:the assumptions are general 2}
Note that Assumptions $\ref{assumption:stabilizable and detectable}$-$\ref{assumption:distance}$ are satisfied by large classes of systems. For example, it was shown in \cite{pasqualetti2012consensus} that Assumption $\ref{assumption:stabilizable and detectable}$ holds if the system dynamics matrix $A$ is row-stochastic and irreducible.\footnote{Note that the matrix $A$ is irreducible if and only if the graph $\mathcal{G}(A)$ is strongly connected \cite{horn1985matrix}.} Assumption $\ref{assumption:distance}$ holds if the system dynamics matrix $A$ is nonnegative and irreducible \cite{horn1985matrix}. Since any row-stochastic matrix is also nonnegative, Assumptions $\ref{assumption:stabilizable and detectable}$-$\ref{assumption:distance}$ hold for any system dynamics matrix $A$ that is row-stochastic and irreducible. Furthermore, using techniques in control theory pertaining to linear structured systems (e.g., \cite{lin1974structural}, \cite{dion2003generic}), one can show that Assumption $\ref{assumption:stabilizable and detectable}$ holds for almost any system dynamics matrix $A$ such that the graph $\mathcal{G}(A)$ is strongly connected, using approaches from \cite{van1991graph,van1999generic}. Specifically, one can consider the system dynamics matrix $A$ to be structured, i.e., each entry of the system dynamics matrix $A$ is either a fixed zero or an independent free parameter (which can attain any real value including zero), where the graph $\mathcal{G}(A)$ is defined according to the free parameters of the structured matrix $A$. One can then show that the set of parameters for which Assumption $\ref{assumption:stabilizable and detectable}$  does not hold has Lebesgue measure zero. Moreover, using similar techniques to those above and the result from \cite{horn1985matrix} that shows that Assumption $\ref{assumption:distance}$ holds for all nonnegative irreducible matrices $A$, one can show that Assumption $\ref{assumption:distance}$ holds for almost any choice of free parameters in the structured matrix $A$ such that the graph $\mathcal{G}(A)$ is strongly connected.  Note that the systems where Assumptions $\ref{assumption:stabilizable and detectable}$-$\ref{assumption:distance}$ hold are not limited to the cases  described above.
\end{remark}
\begin{remark}
We can generalize our analysis to system dynamics matrices $A$ where $\mathcal{G}(A)$ has multiple strongly connected components \cite{bondy1976graph}. Suppose that the input node can only reach (via directed paths in $\mathcal{G}(A)$) nodes that are in the same strongly connected component. Then, under Assumption $\ref{assumption:stabilizable and detectable}$, we only need to consider the strongly connected component of $\mathcal{G}(A)$ that contains the input node, since one can show that the mean square estimation error of the Kalman filter remains zero for the states corresponding to nodes that are not in the strongly connected component containing the input node.
\end{remark}

The first main result of this section is as follows.
\begin{theorem}
Consider a system dynamics matrix $A\in\mathbb{R}^{n\times n}$ with the associated graph $\mathcal{G}(A)=\{\mathcal{X}(A),\mathcal{E}(A)\}$, a single vertex $x_{i_0}\in\mathcal{X}(A)$ that has a stochastic input with variance $\sigma_w^2\in\mathbb{R}_{\ge0}$, the measurement matrix $C=I_n$ (containing all of the individual sensor measurement matrices), and the sensor noise covariance matrix $V=\mathbf{0}_{n\times n}$. Suppose that Assumptions $\ref{assumption:stabilizable and detectable}$-$\ref{assumption:distance}$ hold. For any sensor placement $\mu\in\{0,1\}^n$ such that $\text{supp}(\mu)\neq\emptyset$, denote $\zeta=\mathop{\min}_{j\in\text{supp}(\mu)}l_{i_0j}\ge0$, where $l_{i_0j}$ is the distance from vertex $x_{i_0}$ to vertex $x_j$. The following expressions hold:
\begin{equation}
\Sigma(\mu)=\sigma_{w}^2\displaystyle\sum_{m=0}^{\zeta}A^m B B^T (A^T)^m,
\label{eqn:priori steady state error covariance thm1}
\end{equation}
and
\begin{equation}
\Sigma^*(\mu)=\begin{cases}
&\sigma_{w}^2\displaystyle\sum_{m=0}^{\zeta-1}A^m B B^T (A^T)^m \quad\text{if $\zeta\ge1$},\\
&\mathbf{0} \quad\text{if $\zeta=0$},
\end{cases}
\label{eqn:posteriori steady state error covariance thm1}
\end{equation}
where $\Sigma(\mu)$ (resp., $\Sigma^*(\mu)$) is the steady-state {\it a priori} (resp., {\it a posteriori}) error covariance of the corresponding Kalman filter, and $B=\mathbf{e}_{i_0}$.
\label{thm:steady state error covariance for single input}
\end{theorem}
\begin{IEEEproof}
The existence of $\Sigma(\mu)$ and $\Sigma^*(\mu)$ follows directly from Lemma $\ref{lemma:Anderson optimal filtering}$ and Assumption $\ref{assumption:stabilizable and detectable}$. Considering any sensor placement $\mu$ such that $\zeta\ge1$, i.e., sensors are not placed at the input vertex $x_{i_0}$, we first prove Eq. \eqref{eqn:priori steady state error covariance thm1} by verifying that Eq. \eqref{eqn:priori steady state error covariance thm1} satisfies Eq. \eqref{eqn:DARE}. Note that $C_i=\mathbf{e}^T_i$ for all $x_i\in\mathcal{X}(A)$. Denote $\mathcal{X}_{\mu}\subseteq\mathcal{X}(A)$ as the set of vertices indicated by $\mu$ where sensors are placed and $\mathcal{X}_{\zeta}\subseteq\mathcal{X}(A)$ as the set of vertices that have distance $\zeta$ from the input vertex $x_{i_0}$. Since performing elementary row operations on $C(\mu)$ does not change $\Sigma(\mu)$, we assume without loss of generality that $\mu=\begin{bmatrix}\mu^T_1 & \mu^T_2\end{bmatrix}^T$ such that $\mu_1=\mathbf{1}_{|\mathcal{X}_{\zeta}\cap\mathcal{X}_{\mu}|}$ and $\mu_2\in\{0,1\}^{n-|\mathcal{X}_{\zeta}\cap\mathcal{X}_{\mu}|}$. In other words, $\mu_1$ contains all sensors placed at vertices that have distance $\zeta$ from the input vertex $x_{i_0}$, and $l_{i_0j}>\zeta$ for all $j\in\textrm{supp}(\mu_2)$. The corresponding measurement matrix is given by $C(\mu)=\begin{bmatrix} C(\mu_1)\\ C(\mu_2)\end{bmatrix}$, where $C(\mu_1)\in\mathbb{R}^{|\text{supp}(\mu_1)|\times n}$ and $C(\mu_2)\in\mathbb{R}^{|\text{supp}(\mu_2)|\times n}$. Substituting Eq. \eqref{eqn:priori steady state error covariance thm1} into the right hand side (RHS) of Eq. \eqref{eqn:DARE}, we obtain:
\begin{align}\nonumber
&\text{RHS of Eq. }\eqref{eqn:DARE}\\\nonumber
&=\sigma_w^2\sum_{m=1}^{\zeta+1}A^mBB^T(A^T)^m+\sigma_w^2BB^T-\sigma_w^2A^{\zeta+1}BB^T(A^T)^{\zeta}\\\nonumber
&\qquad\qquad\qquad\quad \times(C(\mu))^T(C(\mu)A^{\zeta}BB^T(A^T)^{\zeta}(C(\mu))^T)^{\dagger}\\
&\qquad\qquad\qquad\qquad\qquad\qquad\times C(\mu)A^{\zeta}BB^T(A^T)^{\zeta+1}\label{eqn:derive priori steady state error cov 1}\\\nonumber
&=\sigma_w^2\sum_{m=0}^{\zeta+1}A^mBB^T(A^T)^m-\sigma_w^2A^{\zeta+1}BB^T(A^T)^{\zeta}\\\nonumber
&\qquad\qquad\times\big[(C(\mu_1))^T\ (C(\mu_2))^T\big] \Big(\begin{bmatrix}C(\mu_1)\\ C(\mu_2)\end{bmatrix}A^{\zeta}BB^T(A^T)^{\zeta}\\\nonumber
&\qquad\quad\times\big[(C(\mu_1))^T \ (C(\mu_2))^T\big]\Big)^{\dagger}\begin{bmatrix}C(\mu_1)\\ C(\mu_2)\end{bmatrix}A^{\zeta}BB^T(A^T)^{\zeta+1}\\\nonumber
&=\sigma_w^2\sum_{m=0}^{\zeta+1}A^mBB^T(A^T)^m -\sigma_w^2A^{\zeta+1}B\\\nonumber
&\qquad\qquad\qquad\qquad\quad\times\big[B^T(A^T)^{\zeta}(C(\mu_1))^T \ \mathbf{0}_{1\times|\text{supp}(\mu_2)|}\big]\\\nonumber
&\qquad\qquad\qquad\quad\times \begin{bmatrix}(C(\mu_1)A^{\zeta}BB^T(A^T)^{\zeta}(C(\mu_1))^T)^{\dagger} & \mathbf{0}\\ \mathbf{0} & \mathbf{0}\end{bmatrix}\\
&\qquad\qquad\qquad\qquad\quad\times\begin{bmatrix}C(\mu_1)A^{\zeta}B\\ \mathbf{0}_{|\text{supp}(\mu_2)|\times 1}\end{bmatrix}B^T(A^T)^{\zeta+1}\label{eqn:derive priori steady state error cov 2}\\\nonumber
&= \sigma_w^2\sum_{m=0}^{\zeta+1}A^mBB^T(A^T)^m-\sigma_w^2A^{\zeta+1}BB^T(A^T)^{\zeta}(C(\mu_1))^T\\
&\times (C(\mu_1)A^\zeta BB^T(A^T)^{\zeta}(C(\mu_1))^T)^{\dagger}C(\mu_1)A^{\zeta}BB^T(A^T)^{\zeta+1},\label{eqn:derive priori steady state error cov 3}
\end{align}
where Eq. \eqref{eqn:derive priori steady state error cov 1} uses the fact that $(A^m)_{ji_0}=0$ for all $j\in\text{supp}(\mu)$ whenever $m\in\{0,1,\dots,\zeta-1\}$, which implies that  $C(\mu)A^mB=\mathbf{0}$ for all $m\in\{0,1,\dots,\zeta-1\}$. Similarly, Eq. \eqref{eqn:derive priori steady state error cov 2} follows from the fact that $C(\mu_2)A^mB=\mathbf{0}$ for all $m\in\{0,1,\dots,\zeta\}$. Denoting $\psi\triangleq C(\mu_1)A^{\zeta}B\in\mathbb{R}^{|\text{supp}(\mu_1)|}$ and noting that $\psi\neq\mathbf{0}$ from Assumption $\ref{assumption:distance}$, one can show that $\psi^T(\psi\psi^T)^{\dagger}\psi=1$. We then have from Eq. \eqref{eqn:derive priori steady state error cov 3}:
\begin{align*}
&\text{RHS of Eq. }\eqref{eqn:DARE}\\
&=\sigma_w^2\sum_{m=0}^{\zeta+1}A^mBB^T(A^T)^m-\sigma_w^2A^{\zeta+1}BB^T(A^T)^{\zeta+1}\\
&=\sigma_w^2\sum_{m=0}^{\zeta}A^mBB^T(A^T)^m.
\end{align*}
Since $\sigma_w^2\sum_{m=0}^{\zeta}A^mBB^T(A^T)^m\succeq\mathbf{0}$, we know from Lemma $\ref{lemma:Anderson optimal filtering}$ that the limit $\Sigma(\mu)$ is given by Eq. \eqref{eqn:priori steady state error covariance thm1}. We then obtain from Eq. \eqref{eqn:coupled DARE} that the limit $\Sigma^*(\mu)$ is given by Eq. \eqref{eqn:posteriori steady state error covariance thm1} (when $\zeta\ge1$).

Next, we consider any sensor placement $\mu$ such that $\zeta=0$, i.e., a sensor is placed at the input vertex $x_{i_0}$. Using similar arguments to those above, we can also show that Eq. \eqref{eqn:priori steady state error covariance thm1}-\eqref{eqn:posteriori steady state error covariance thm1} hold when $\zeta=0$. This completes the proof of the theorem.\end{IEEEproof}

\begin{figure}[hp]
\vspace{-0.3cm}
\centering
\captionsetup{justification=centering}
\includegraphics[width=0.6\linewidth]{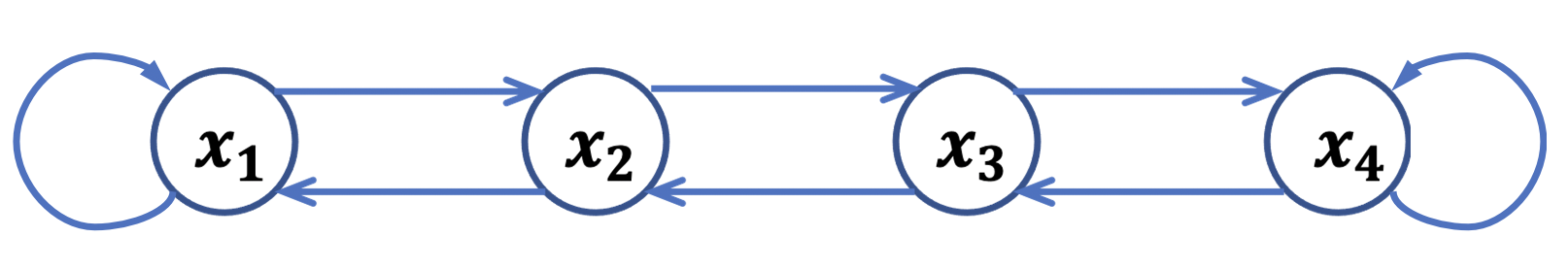}
\center
\vspace{-0.4cm}
\caption{Graph for Example $1$}
\label{fig:thm1 example}
\vspace{-0.6cm}
\end{figure}

To verify the results in Theorem $\ref{thm:steady state error covariance for single input}$, let us consider the following example.
\begin{example}
Consider the graph in Fig. $\ref{fig:thm1 example}$, where $x_2$ is the input node (i.e., $B=\mathbf{e}_2$) with variance $\sigma_w^2=1$.  Suppose $A=\left[\begin{smallmatrix}0.5 & 2.1 & 0 & 0\\ 0.3 & 0 & 1.5 & 0\\ 0 & 0.6 & 0 & 0.5\\ 0 & 0 & -0.8 & 1\end{smallmatrix}\right]$, $C=I_4$ and $V=\mathbf{0}_{4\times 4}$. Denote $\mu_2=[0\ 1\ 0\ 0]^T$ and $\mu_4=[0\ 0\ 0\ 1]^T$. It can be verified that $\Sigma(\mu_2)=\left[\begin{smallmatrix}0 & 0 & 0 & 0\\0 & 1 & 0 & 0\\0 & 0 & 0 & 0\\0 & 0 & 0 & 0\end{smallmatrix}\right]=BB^T$, $\Sigma^*(\mu_2)=\mathbf{0}_{4\times4}$, $\Sigma(\mu_4)=\left[\begin{smallmatrix}5.5125 & 1.6065 & 1.26 & -0.504\\1.6065 & 3.3409 &  0 & -0.7344\\1.26 &  0 & 0.36 &  0\\-0.504 & -0.7344 &  0 &  0.2304\end{smallmatrix}\right]=\sum_{m=0}^{2}A^mBB^T(A^T)^m$ and $\Sigma^*(\mu_4)=\left[\begin{smallmatrix}4.41 & 0 & 1.26 & 0\\0 & 1 & 0 & 0\\1.26 & 0 & 0.36 & 0\\0 & 0 & 0 & 0\end{smallmatrix}\right]=\sum_{m=0}^{1}A^mBB^T(A^T)^m$, as provided by Theorem $\ref{thm:steady state error covariance for single input}$.
\end{example}

\subsection{An Optimal Solution to GKFSP}\label{sec:optimal solution to GKFSP}
Using the above discussions, we give the following result that characterizes an optimal solution to GKFSP (Problem $\ref{pro:GKFSP problem}$).
\begin{theorem}
\label{thm:optimal solution to GKFSP}
Supposing that Assumptions $\ref{assumption:stabilizable and detectable}$-$\ref{assumption:distance}$ hold, an optimal solution, denoted as $\mu^*$, to the priori (resp., posteriori) GKFSP problem is to place a  single sensor at a vertex $x_j$ in order to minimize $l_{i_0j}$, i.e., the distance from the input vertex $x_{i_0}$ to $x_j$, while satisfying the budget constraint.
\end{theorem}
\begin{IEEEproof}
Under Assumptions $\ref{assumption:stabilizable and detectable}$-$\ref{assumption:distance}$, we first note from Eq. \eqref{eqn:priori steady state error covariance thm1}-\eqref{eqn:posteriori steady state error covariance thm1} that the {\it a priori} and the {\it a posteriori} error covariance matrices only depend on $\zeta$, i.e., the shortest distance from the input node to the sensor nodes. Hence, it is sufficient to consider sensor placements $\mu\in\{0,1\}^{n}$ such that $|\text{supp}(\mu)|=1$ in terms of minimizing the trace of the {\it a priori} (resp., {\it a posteriori}) steady-state error covariance of the Kalman filter. Moreover, we know from Eq. \eqref{eqn:priori steady state error covariance thm1} in Theorem $\ref{thm:steady state error covariance for single input}$ that $\Sigma(\mu)=\sigma_w^2\sum_{m=0}^{\zeta}A^m B B^T (A^T)^m$, where $\zeta=\mathop{\min}_{j\in\text{supp}(\mu)}l_{i_0j}$. Noting that the matrix $A^mBB^T(A^T)^m$ is positive semi-definite for all $m\in\mathbb{Z}_{\ge0}$, it follows that $\text{trace}(A^mBB^T(A^T)^m)\ge0$ for all $m\in\mathbb{Z}_{\ge0}$. Hence, $\text{trace}(\Sigma(\mu))$ is minimized by finding a sensor placement $\mu^*$ with $|\text{supp}(\mu^*)|=1$ such that $\zeta$ is minimized while satisfying the budget constraint. Using similar arguments, we can show that $\mu^*$ is also an optimal solution to the posteriori GKFSP problem.
\end{IEEEproof}

Based on Theorem $\ref{thm:optimal solution to GKFSP}$, we can find the optimal solution $\mu^*$ to the priori (resp., posteriori) GKFSP problem using polynomial-time algorithms such as the Breadth-First Search (BFS) algorithm which runs in time $O(n+|\mathcal{E}(A)|)$ \cite{cormen2009introduction}. 

\subsection{An Optimal Solution to GKFSA}\label{sec:solution to GKFSA}
Given a sensor placement $\mu$, we know from the insights obtained above for GKFSP that the steady-state {\it a priori} and the {\it a posteriori} error covariances of the Kalman filter (after an attack that removes some of those sensors) only depend on the surviving sensors that have the shortest distance from the input vertex $x_{i_0}$. We then have the following result whose proof is similar to that of Theorem $\ref{thm:optimal solution to GKFSP}$ and is thus omitted.
\begin{theorem}
\label{thm:optimal solution to GKFSA}
Suppose that Assumptions $\ref{assumption:stabilizable and detectable}$-$\ref{assumption:distance}$ hold. Given a sensor placement $\mu$, an optimal solution, denoted as $\nu^*$, to the priori (resp., posteriori) GKFSA problem can be found by maximizing the shortest distance from the input vertex $x_{i_0}$ to the surviving sensors, i.e., solving the following optimization problem
\begin{equation}
\label{eqn:equivalent shortest path problem attack}
\begin{split}
&\mathop{\max}_{\nu\in\{0,1\}^n}\mathop{\min}_{j\in\text{supp}(\mu\setminus\nu)}l_{i_0j}\\
&s.t.\ f^T\nu\le F,
\end{split}
\end{equation}
where $l_{i_0j}$ is the distance from vertex $x_{i_0}$ to vertex $x_j$, and $l_{i_0j}=+\infty$ if $\text{supp}(\mu\setminus\nu)=\emptyset$.
\end{theorem}

The optimal solution $\nu^*$ to the priori (resp., posteriori) GKFSA problem described by Theorem $\ref{thm:optimal solution to GKFSA}$ can be found as follows. Given a sensor placement $\mu$, the adversary starts by inspecting the placed sensors (indicated by $\mu$) that have the shortest distance from the input vertex $x_{i_0}$. The adversary will remove all of these sensors if the sum of the corresponding sensor attack costs is less than or equal to the budget constraint $F$, and terminate the process if otherwise. The above process is then repeated for the placed sensors that have the second shortest distance from the input vertex $x_{i_0}$, based on the remaining budget. This process continues with the placed sensors that have the third shortest distance from the input vertex $x_{i_0}$, etc.  Hence, polynomial-time algorithms such as the BFS algorithm can be used to find the optimal sensor attack $\nu^*$ for the adversary in time $O(n+|\mathcal{E}(A)|)$. 

\section{Solving the RGKFSP problem}\label{sec:alg for RGKFSP}
We now turn to the RGKFSP problem (Problem $\ref{pro:RGKFSP problem}$). Recall that Theorem $\ref{thm:optimal solution to GKFSP}$ showed that it is enough to consider only sensor placements $\mu$ with $|\text{supp}(\mu)|=1$ for the GKFSP problem (i.e., the system designer does not necessarily need to utilize all of the sensor placement budget $H$). However, an optimal sensor placement $\mu^*$ for the RGKFSP problem does not necessarily satisfy $|\text{supp}(\mu^*)|=1$, since the adversary could have enough budget to remove the single sensor placed by the system designer, which causes the trace of the {\it a priori} (resp., {\it a posteriori}) error covariance of the Kalman filter to be infinite (if the system dynamics matrix $A$ is not stable). Note that the steady-state {\it a priori} and the {\it a posteriori} error covariance matrices of the Kalman filter (after the attack) only depend on the surviving sensors that have the shortest distance from the input vertex $x_{i_0}$. Using similar arguments to those for Theorems $\ref{thm:optimal solution to GKFSP}$-$\ref{thm:optimal solution to GKFSA}$, we have that an optimal solution to the RGKFSP problem can be found by minimizing the shortest distance from the input vertex $x_{i_0}$ to the sensors after the corresponding optimal sensor attack, and a sensor placement $\mu^*$ is optimal for the priori RGKFSP problem if and only if it is optimal for the posteriori RGKFSP problem.

We thus focus on the priori RGKFSP problem in this section. Although we provided polynomial-time algorithms to solve the GKFSP and GKFSA problems, we will show that the RGKFSP problem is NP-hard, i.e., there exist classes of the RGKFSP problem  that cannot be solved by any polynomial-time algorithm if P$\neq$NP. To do this, we first recall from Remark $\ref{remark:the assumptions are general 2}$ that Assumptions $\ref{assumption:stabilizable and detectable}$-$\ref{assumption:distance}$ hold for any system dynamics matrix $A$ that is row-stochastic and irreducible. Therefore, Eq. \eqref{eqn:priori steady state error covariance thm1} and Eq. \eqref{eqn:posteriori steady state error covariance thm1} in Theorem $\ref{thm:steady state error covariance for single input}$ also hold for such $A$ matrices.

To show that the RGKFSP problem is NP-hard, we reduce the subset sum problem \cite{garey1979computers} to RGKFSP.
\begin{definition}
\label{def:subset sum instance}
An instance of the subset sum problem is given by a finite set $U$ and a positive integer $K$, where each $s\in U$ has a size $\kappa(s)\in\mathbb{Z}_{>0}$. 
\end{definition}

We use the following result from \cite{garey1979computers}.
\begin{lemma}
\label{lemma:subset sum np-complete}
Given an instance of the subset sum problem as described in Definition $\ref{def:subset sum instance}$, the problem of determining whether there is a subset $U'\subseteq U$ such that $\sum_{s\in U'}\kappa(s)=K$ is NP-complete. 
\end{lemma}

We are now in place to prove the following result.
\begin{theorem}
\label{thm:RGKFSP np-hard subset sum}
The RGKFSP problem is NP-hard even when both of the following two conditions are satisfied: $(1)$ the sensor placement cost and the sensor attack cost satisfy $h_i=f_i$ for all $i\in\{1,2,\dots,n\}$, and $(2)$ there is a feasible sensor placement for the system designer.
\end{theorem}
\begin{IEEEproof}
We prove the result by giving a polynomial-time reduction from the subset sum problem. Consider any instance of the subset sum problem  defined in Definition $\ref{def:subset sum instance}$. Denote $U=\{s_1,s_2,\dots,s_{|U|}\}$. Denote the number of bits of the binary representation of the positive integer $K$ as $b(K)$, i.e., $b(K)\triangleq\lfloor\text{log}_2(K)\rfloor+1$. We then construct an instance of the priori RGKFSP problem as follows. The system dynamics matrix $A\in\mathbb{R}^{(|U|+b(K))\times(|U|+b(K))}$ is chosen such that the graph $\mathcal{G}(A)$ is an undirected path of length $|U|+b(K)-1$. Specifically, we set $A_{ij}=A_{ji}=\frac{1}{3}$ for all $i\in\{1,2,\dots,|U|+b(K)-1\}$ and $j=i+1$, $A_{mm}=\frac{1}{3}$ for all $m\in\{2,3,\dots,|U|+b(K)-1\}$, $A_{mm}=\frac{2}{3}$ for all $m\in\{1,|U|+b(K)\}$, and all the other entries in $A$ are zero. The vertex $x_1$ is set as the only vertex that has the stochastic input with variance $\sigma_w^2=1$. The sensor placement cost vector is set as $h_i=\kappa(s_i)$ for all $i\in\{1,2,\dots,|U|\}$, and $h_i=2^{i-|U|-1}$ for all $i\in\{|U|+1,|U|+2,\dots,|U|+b(K)\}$. The sensor attack cost is set as $f_i=h_i$ for all $i\in\{1,2,\dots,|U|+b(K)\}$. Note that the sensor placement vector and the sensor attack vector are given by $\mu\in\{0,1\}^{|U|+b(K)}$ and $\nu\in\{0,1\}^{|U|+b(K)}$, respectively. The sensor placement budget of the system designer is set as $H=K$, and the sensor attack budget of the adversary is set as $F=K-1$. We also note that the matrix $A$ that we constructed is row-stochastic and irreducible. Therefore, Eq. \eqref{eqn:priori steady state error covariance thm1} in Theorem $\ref{thm:steady state error covariance for single input}$ holds for the $A$ matrix that we constructed. We claim that the answer to the given subset sum instance is ``yes'' if and only if the optimal solution to the constructed instance of the priori RGKFSP problem, denoted as $\mu^*$, satisfies $\text{trace}(\Sigma(\mu^*\setminus\nu^*))\le\text{trace}(\sum_{i=0}^{|U|-1}A^iBB^TA^i)$,  where $\nu^*$ is the optimal sensor attack given $\mu^*$.

Suppose that the answer to the given subset sum instance is ``yes'', i.e., there exists $U'\subseteq U$ such that $\sum_{s\in U'}\kappa(s)=K$. It follows that for the instance of the priori RGKFSP problem as constructed above, there exists a sensor placement vector $\tilde{\mu}$ such that $\sum_{i=1}^{|U|}h_i \tilde{\mu}_i=K\le H$. Therefore, for any sensor attack $\tilde{\nu}$ that satisfies the sensor attack budget constraint, i.e., $\sum_{i=1}^{|U|}f_i \tilde{\nu}_i\le F=K-1$, we have $\text{supp}(\tilde{\mu}\setminus\tilde{\nu})\cap\{1,\dots,|U|\}\neq\emptyset$, which implies that there exists $j\in\{1,\dots,|U|\}$ such that $j\in\text{supp}(\tilde{\mu}\setminus\tilde{\nu})$. Noting that $A^mBB^TA^m\succeq\mathbf{0}$ for all $m\in\mathbb{Z}_{\ge0}$ and $l_{1j}=j-1\le l_{1|U|}=|U|-1$, it then follows from Eq. \eqref{eqn:priori steady state error covariance thm1} that $\text{trace}(\Sigma(\tilde{\mu}\setminus\tilde{\nu}))\le\text{trace}(\sum_{i=0}^{j-1}A^iBB^TA^i)\le\text{trace}(\sum_{i=0}^{|U|-1}A^iBB^TA^i)$, for any sensor attack $\tilde{\nu}$ such that $\sum_{i=1}^{|U|}f_i \tilde{\nu}_i\le F$. Since $\text{trace}(\Sigma(\mu^*\setminus\nu^*))\le\text{trace}(\Sigma(\tilde{\mu}\setminus\tilde{\nu}))$, we have $\text{trace}(\Sigma(\mu^*\setminus\nu^*))\le\text{trace}(\sum_{i=0}^{|U|-1}A^iBB^TA^i)$. 

Conversely, suppose that the answer to the subset sum instance is ``no'', i.e., for any $U'\subseteq U$, we have $\sum_{s\in U'}\kappa(s)\neq K$. Considering the instance of the priori RGKFSP problem we constructed, for any sensor placement vector $\mu$ such that $\sum_{i=1}^{|U|+b(K)}h_i \mu_i\le H=K$, we have $\sum_{i=1}^{|U|}h_i \mu_i\neq K$,  which implies $\sum_{i=1}^{|U|}h_i \mu_i\le K-1$. Denote $\sum_{i=1}^{|U|}h_i \mu_i\triangleq K_{|U|}$. Therefore, for any sensor placement vector $\mu$ with $\sum_{i=1}^{|U|+b(K)}h_i\mu_i\le H$, there exists an attack $\hat{\nu}$ such that $\sum_{i=1}^{|U|}f_i \hat{\nu}_i=K_{|U|}\le K-1$, which implies $\text{supp}(\mu\setminus\hat{\nu})\cap\{1,2,\dots,|U|\}=\emptyset$. Moreover, note that $K=H\ge H-K_{|U|}>F-K_{|U|}$. Since we set the sensor placement cost vector and the sensor attack cost vector to satisfy $h_i=f_i=2^{i-|U|-1}$ for all $i\in\{|U|+1,|U|+2,\dots,|U|+b(K)\}$, where $b(K)$ is the number of bits for the binary representation of $K$, we have that for any $U'\subseteq U$, there exists $\bar{U}'\subseteq\{|U|+1,|U|+2,\dots,|U|+b(K)\}$ such that $\sum_{s\in U'}\kappa(s)+\sum_{i\in \bar{U}'}h_i=H$. Therefore, the system designer can always use all the sensor placement budget by placing sensors at an appropriate subset of the vertices in the vertex set $\{x_{|U|+1},x_{|U|+2},\dots,x_{|U|+b(K)}\}$ and guarantee to have at least one sensor left after any attack that satisfies the sensor attack budget constraint. Formally, we have that for any sensor placement $\mu$ with $\sum_{i=1}^{|U|+b(K)}h_i\mu_i=H$, there exists $j'\in\{|U|+1,\dots,|U|+b(K)\}$ such that $j'\in\text{supp}(\mu\setminus\nu)$, where $\nu$ is any sensor attack satisfying the sensor attack budget constraint. Meanwhile, any sensor placement $\mu$ such that $\sum_{i=1}^{|U|+b(K)}h_i\mu_i< H$ is not a feasible sensor placement. Therefore, there is always a feasible sensor placement for the system designer under the constructed instance of the priori RGKFSP problem when the answer to the subset sum instance is ``no". Note that the matrix $A^mBB^TA^m\succeq\mathbf{0}$ for all $m\in\mathbb{Z}_{\ge0}$ and $l_{1j'}=j'-1\ge l_{1|U|}+1= |U|$. Combining the arguments above together, it then follows from Eq. \eqref{eqn:priori steady state error covariance thm1} that for any $\mu$ such that $\sum_{i=1}^{|U|+b(K)}h_i \mu_i=H$, we have $\text{trace}(\Sigma(\mu\setminus\nu))\ge\text{trace}(\sum_{i=0}^{j'-1}A^iBB^TA^i)\ge\text{trace}(\sum_{i=0}^{|U|}A^iBB^TA^i)$, where $\nu$ is any sensor attack satisfying the sensor attack budget constraint. Since $(A^m)_{11}>0$ for all $m\in\mathbb{Z}_{\ge0}$, we have $\text{trace}(A^{|U|}BB^TA^{|U|})>0$ and thus $\text{trace}(\Sigma(\mu\setminus\nu))>\text{trace}(\sum_{i=0}^{|U|-1}A^iBB^TA^i)$. Since the above arguments hold for any $\mu$ with $\sum_{i=1}^{|U|+b(K)}h_i\mu_i=H$, they also hold for the optimal solution $\mu^*$ to the constructed priori RGKFSP instance, i.e., $(\Sigma(\mu^*\setminus\nu^*))>\text{trace}(\sum_{i=0}^{|U|-1}A^iBB^TA^i)$, where $\nu^*$ is the optimal sensor attack given $\mu^*$. This completes the proof of the claim above.

Since the subset sum problem is NP-complete and RGKFSP $\notin$ NP, we conclude that RGKFSP is NP-hard even under the additional conditions as stated.
\end{IEEEproof}

\subsection{An Algorithm for RGKFSP}
It follows directly from Theorem $\ref{thm:RGKFSP np-hard subset sum}$ that there is no polynomial-time algorithm that would solve all instances of RGKFSP if P$\neq$NP. However, we now provide a pseudo-polynomial-time algorithm\footnote{A pseudo-polynomial-time algorithm is an algorithm that runs in time that is bounded by a polynomial in the largest integer in its input \cite{garey1979computers}.} (Algorithm $\ref{algorithm:DP for RGKFSP}$) for RGKFSP by relating it to the knapsack problem defined as follows.
\begin{definition}
\label{def:knapsack}
Given a finite set $U\triangleq\{s_1,s_2,\dots,s_{|U|}\}$, a size $\kappa(s_i)\in\mathbb{Z}_{>0}$ and a value $\phi(s_i)\in\mathbb{Z}_{>0}$ for each $i\in\{1,2,\dots,|U|\}$, and a positive integer $K$, the knapsack problem is to find an indicator vector $\pi\in\{0,1\}^{|U|}$ that solves 
\begin{equation}
\label{eqn:knapsack}
\begin{split}
\begin{split}
&\mathop{\max}_{\pi\in\{0,1\}^{|U|}}\sum_{i=1}^{|U|}\phi(s_i)\pi_i\\
&s.t.\ \sum_{i=1}^{|U|}\kappa(s_i)\pi_i\le K.
\end{split}
\end{split}
\end{equation}
Denoting the given instance of knapsack as a tuple $\{\phi,\kappa,K\}$,where $\phi\triangleq(\phi(s_1), \phi(s_2),\dots, \phi(s_{|U|}))$ and $\kappa\triangleq(\kappa(s_1),\kappa(s_2),\dots, \kappa(s_{|U|}))$,\footnote{Note that the elements in $\phi$ and $\kappa$ are ordered, and the $i$th element of $\phi$ (resp., $\kappa$) corresponds to the value (resp., weight) of $s_i\in U$ for all $i\in\{1,\dots,|U|\}$. The dependency of $\{\phi,\kappa,K\}$ on $U$ is dropped since each element of $\phi$ (resp., $\kappa$) represents an element of $U$.} the corresponding optimal indicator vector for \eqref{eqn:knapsack} is denoted as $\pi^*(\phi,\kappa,K)$, and the corresponding optimal value of the objective function in \eqref{eqn:knapsack} is denoted as $z(\phi,\kappa,K)$.
\end{definition}

The steps of Algorithm $\ref{algorithm:DP for RGKFSP}$ for RGKFSP are as follows. Algorithm $\ref{algorithm:DP for RGKFSP}$ starts by relabeling the input vertex as vertex $x_1$ and relabeling the other vertices in terms of a non-decreasing order of the distances from the vertex $x_1$ (Lines $1$-$2$). Denoting $l_{\text{max}}\triangleq\mathop{\max}_{x_j\in\mathcal{X}(A)}l_{i_0j}$, Algorithm $\ref{algorithm:DP for RGKFSP}$ then finds the smallest $m\in\{0,1,\dots,l_{\text{max}}\}$ such that by placing sensors (under the budget constraint) solely at nodes that have distances less than or equal to $m$ from $x_1$ (after the relabeling), the sum of the sensor attack costs of the placed sensors is greater than the sensor attack budget, i.e., there is at least one sensor that survives the corresponding optimal sensor attack. This is done by iteratively solving a knapsack problem at increasingly longer distances from the input node, where at each distance, the goal is to find a set of sensor locations that fits within the sensor placement budget constraint $H$ but maximizes the sum of the sensor attack costs.  Algorithm $\ref{algorithm:DP for RGKFSP}$ returns $\mu=\mathbf{0}_{n\times 1}$ if there is no feasible sensor placement. We now prove that Algorithm $\ref{algorithm:DP for RGKFSP}$ returns an optimal solution to RGKFSP.

\begin{algorithm}
\textbf{Input:} An instance of the RGKFSP problem.\\
\textbf{Output:} A sensor placement $\mu\in\{0,1\}^n$.
\caption{Algorithm for RGKFSP}\label{algorithm:DP for RGKFSP}
\begin{algorithmic}[1]
\State Find the distance $l_{i_0j}$ for all $x_j\in\mathcal{X}(A)\setminus\{x_{i_0}\}$ via BFS and denote $l_{\text{max}}\triangleq\mathop{\max}_{x_j\in\mathcal{X}(A)}l_{i_0j}$.
\State Relabel the vertices of $\mathcal{G}(A)$ such that $x_1$ is the input vertex and $l_{1j}\le l_{1t}$ for all $x_j,x_t\in\mathcal{X}(A)\setminus\{x_1\}$ with $j\le t$.
\State $\mu=\mathbf{0}_{n\times 1}$
\For{$m=0$ to $l_{\text{max}}$}
\State Find $j_m\triangleq\mathop{\max}\{j:l_{1j}= m,x_j\in\mathcal{X}(A)\}$.
\State Find $\pi^*\big((f_1,\dots,f_{j_m}),(h_1,\dots,h_{j_m}), H\big)$
\If {$z\big((f_1,\dots,f_{j_m}),(h_1,\dots,h_{j_m}), H\big)>F$}
\State $[\mu_1 \ \cdots  \ \mu_{j_m}]^T=$
\Statex \qquad\quad$\pi^*\big( (f_1,\dots,f_{j_m}),(h_1,\dots,h_{j_m}), H\big)$
\EndIf
\EndFor
\State\qquad\quad\Return {$\mu$}
\end{algorithmic}
\end{algorithm}

\begin{theorem}
\label{thm:optimal of dynamic programming for RGKFSP}
Under Assumptions $\ref{assumption:stabilizable and detectable}$-$\ref{assumption:distance}$, Algorithm $\ref{algorithm:DP for RGKFSP}$ returns an optimal solution to the RGKFSP problem.
\end{theorem}
\begin{IEEEproof}
Denote an optimal solution to the RGKFSP problem as $\mu^*$ and denote the solution returned by Algorithm $\ref{algorithm:DP for RGKFSP}$ as $\mu'$. Suppose that $\mu'$ is a feasible sensor placement. Suppose that the vertices in $\mathcal{G}(A)$ are  relabeled as indicated by Lines $1$-$2$ in Algorithm $\ref{algorithm:DP for RGKFSP}$, i.e., vertex $x_1$ is labeled as the input vertex and the other vertices are labeled in terms of a non-decreasing order of the distances from vertex $x_1$ (note that the relabeling of the vertices does not change the optimal solution to the RGKFSP problem other than permuting it). Assume for the sake of contradiction that $\text{trace}(\Sigma(\mu^*\setminus\nu^*))<\text{trace}(\Sigma(\mu'\setminus\nu'))$, where $\nu^*$ and $\nu'$ are the optimal sensor attacks given $\mu^*$ and $\mu'$, respectively. Denote $j^*\triangleq\mathop{\max}J$ and $j'\triangleq\mathop{\max}J'$, where $J\triangleq\mathop{\arg}{\min}_{m\in\text{supp}(\mu^*\setminus\nu^*)}l_{1m}$ and $J'\triangleq\mathop{\arg}{\min}_{m\in\text{supp}(\mu'\setminus\nu')}l_{1m}$. In other words, among those sensors that are closest to the input vertex in $\text{supp}(\mu^*\setminus\nu^*)$ (resp., $\text{supp}(\mu'\setminus\nu')$), $j^*$ (resp., $j'$) is the largest index.  Noting that $\sum_{m=1}^{j^*}f_m\mu^*_m>F$ (otherwise the optimal sensor attack $\nu^*$ given $\mu^*$ would remove the sensor placed at vertex $x_{j^*}$ as argued previously in Section $\ref{sec:solution to GKFSA}$), it follows from Definition $\ref{def:knapsack}$ that $z\big((f_1,\dots,f_{j^*}),(h_1,\dots,h_{j^*}),H\big)>F$, which implies that $z\big((f_1,\dots,f_{j_m}),(h_1,\dots,h_{j_m}),H\big)>F$, where $j_m$ is defined in Line $5$ of Algorithm $\ref{algorithm:DP for RGKFSP}$ with $m=l_{1j^*}$. We then know from the definition of Algorithm $\ref{algorithm:DP for RGKFSP}$ that the sensor placement $\mu'$ returned by Algorithm $\ref{algorithm:DP for RGKFSP}$ would satisfy $j'\le j_m$, which implies that $l_{1j'}\le l_{1j^*}$ (by the way that Algorithm $\ref{algorithm:DP for RGKFSP}$ relabels the vertices). Moreover, we have from Theorem $\ref{thm:steady state error covariance for single input}$ that
\begin{equation}
\label{eqn:trace of surviving sensors 1}
\Sigma(\mu^*\setminus\nu^*)=\sigma_w^2\displaystyle\sum_{m=0}^{l_{1j^*}}A^m B B^T (A^T)^m,
\end{equation}
and
\begin{equation}
\label{eqn:trace of surviving sensors 2}
\Sigma(\mu'\setminus\nu')=\sigma_w^2\displaystyle\sum_{m=0}^{l_{1j'}}A^m B B^T (A^T)^m,
\end{equation}
hold  under Assumptions $\ref{assumption:stabilizable and detectable}$-$\ref{assumption:distance}$. Since the matrix $A^mBB^T(A^T)^m\succeq\mathbf{0}$ for all $m\in\mathbb{Z}_{\ge0}$, we have from the assumption $\text{trace}(\Sigma(\mu^*\setminus\nu^*))<\text{trace}(\Sigma(\mu'\setminus\nu'))$ and Eq. \eqref{eqn:trace of surviving sensors 1}-\eqref{eqn:trace of surviving sensors 2} that $l_{1j^*}<l_{1j'}$. Thus, we get a contradiction. 

We then suppose that the solution $\mu'$ returned by Algorithm $\ref{algorithm:DP for RGKFSP}$ is not feasible, i.e., $\text{supp}(\mu'\setminus\nu')=\emptyset$. Again, we assume that $\text{trace}(\Sigma(\mu^*\setminus\nu^*))<\text{trace}(\Sigma(\mu'\setminus\nu'))$, i.e., $\text{supp}(\mu^*\setminus\nu^*)\neq\emptyset$. Via similar arguments to those above, we have that there exists $j^{*\prime}\in\{1,\dots,n\}$ such that $z\big((f_1,\dots,f_{j^{*\prime}}),(h_1,\dots,h_{j^{*\prime}}),H\big)>F$, which implies that $z\big((f_1,\dots,f_{j_m}),(h_1,\dots,h_{j_m}),H\big)>F$, where $j_m$ is defined in Line $5$ of Algorithm $\ref{algorithm:DP for RGKFSP}$ with $m=l_{1j^{*\prime}}$. Therefore, Algorithm $\ref{algorithm:DP for RGKFSP}$ would also return a solution $\mu'$ such that $\text{supp}(\mu'\setminus\nu')\neq\emptyset$, which is a contradiction. We then conclude that $\text{trace}(\Sigma(\mu^*\setminus\nu^*))=\text{trace}(\Sigma(\mu'\setminus\nu'))$, i.e., Algorithm $\ref{algorithm:DP for RGKFSP}$ returns an optimal solution to the RGKFSP problem. 
\end{IEEEproof}

Since the knapsack problem is NP-hard, there is no polynomial-time algorithm to solve it optimally (if P$\neq$NP) \cite{garey1979computers}. Various algorithms exist to approximate or optimally solve it, including greedy algorithms, linear programming relaxation and dynamic programming \cite{kellerer2004knapsack}. When implementing Algorithm $\ref{algorithm:DP for RGKFSP}$, we can use existing algorithms for knapsack to find $\pi^*\big((f_1,\dots,f_{j_m}),(h_1,\dots,h_{j_m}), H\big)$ in Line $6$ and $z\big((f_1,\dots,f_{j_m}),(h_1,\dots,h_{j_m}), H\big)$ in Line $7$ when we range $m$ from $0$ to $l_{\text{max}}$. Specifically, we call a pseudo-polynomial-time algorithm for knapsack (that solves it optimally) at most $l_{\text{max}}+1$ times to achieve this. For example, a typical dynamic programming approach for knapsack finds $\pi^*\big((f_1,\dots,f_{j_m}),(h_1,\dots,h_{j_m}), H\big)$ and $z\big((f_1,\dots,f_{j_m}),(h_1,\dots,h_{j_m}), H\big)$ in time $O(j_mH)$ for each $m\in\{0,\dots,l_{\text{max}}\}$ \cite{kellerer2004knapsack}. Since BFS runs in time $O(n+|\mathcal{E}(A)|)$, Algorithm $\ref{algorithm:DP for RGKFSP}$ runs in time $O(l_{\text{max}}nH+n+|\mathcal{E}(A)|)$.

\section{Noisy Sensor Measurement Case}\label{sec:noisy measurement}
The results we obtained so far hold under the assumption that $V=\mathbf{0}_{n\times n}$. In this section, we provide a bound on the suboptimality of the proposed strategies when there is sensor measurement noise. We will use the following result whose proof is in the appendix.

\begin{lemma}
\label{lemma:small noise case}
Consider a system dynamics matrix $A\in\mathbb{R}^{n\times n}$, an input matrix $B\in\mathbb{R}^{n\times n_1}$, a sensor measurement matrix $C\in\mathbb{R}^{n_2\times n}$, an input covariance matrix $W\in\mathbb{S}^{n_1}_{+}$, and a sensor measurement noise covariance matrix $\tilde{V}\in\mathbb{S}^{n_2}_{+}$. Suppose that the pair $(A,BW^{1/2})$ is stabilizable and the pair $(A,C)$ is detectable. Denote $\tilde{\Sigma}$ (resp., $\tilde{\Sigma}^*$) as the steady-state {\it a priori} (resp., {\it a posteriori}) error covariance of the Kalman filter corresponding to the measurement noise covariance $\tilde{V}$, and denote $\Sigma$ (resp., $\Sigma^*$) as the steady-state {\it a priori} (resp., {\it a posteriori}) error covariance of the corresponding Kalman filter when $V=\mathbf{0}_{n_2\times n_2}$. Then, 
$\tilde{\Sigma}\preceq\Sigma+E$
and
$\tilde{\Sigma}^*\preceq\Sigma^*+(I_n-LC)E$,
where $E$ is given by
\begin{equation}
\label{eqn:extra term on priori upper bound}
E\triangleq\sum_{m=0}^{\infty}(A-KC)^mK\tilde{V}K^T((A-KC)^T)^m,
\end{equation}
with
$K\triangleq A\Sigma C^T(C\Sigma C^T)^{-1}$ and 
$L\triangleq\Sigma C^T(C\Sigma C^T)^{-1}$.\footnote{The inverses are interpreted as the Moore-Penrose pseudo-inverses  if the arguments are not invertible \cite{anderson1979optimal}.}
\end{lemma}

Note that $E$ exists and is finite since the matrix $A-KC$ is stable. See the proof in the appendix for more details. We have the following result for the GKFSP problem. 
\begin{theorem}
\label{thm:noisy measurement case}
Suppose that Assumptions $\ref{assumption:stabilizable and detectable}$-$\ref{assumption:distance}$ hold. Let $\tilde{\Sigma}(\mu)$ (resp., $\tilde{\Sigma}^*(\mu)$) be the steady-state {\it a priori} (resp., {\it a posteriori}) error covariance matrix of the Kalman filter associated with $\mu$ when $V=\tilde{V}\in\mathbb{S}^n_{+}$. Denote $\tilde{\mu}^*_1$ (resp., $\tilde{\mu}^*_2$) as the optimal solution to the priori (resp., posteriori) GKFSP problem when $V=\tilde{V}$, and denote $\mu^*$ as the optimal solution to the priori (resp., posteriori) GKFSP problem when $V=\mathbf{0}_{n\times n}$. Then, $\text{trace}(\tilde{\Sigma}(\mu^*))\le\text{trace}(\tilde{\Sigma}(\tilde{\mu}_1^*))+\text{trace}(E(\mu^*))$ and $\text{trace}(\tilde{\Sigma}^*(\mu^*))\le\text{trace}(\tilde{\Sigma}^*(\tilde{\mu}_2^*))+\text{trace}((E^*(\mu^*))$,  where $E(\mu^*)$ and $L(\mu^*)$ are defined in Lemma $\ref{lemma:small noise case}$ with $C=C(\mu^*)$, and $E^*(\mu^*)\triangleq(I_n-L(\mu^*)C(\mu^*))E(\mu^*)$.
\end{theorem}
\begin{IEEEproof}
First, we know from Lemma $\ref{lemma:small noise case}$ that $\tilde{\Sigma}(\mu^*)\preceq\Sigma(\mu^*)+E(\mu^*)$, where $\Sigma(\mu^*)$ is the steady-state {\it a priori} error covariance of the Kalman filter corresponding to $\mu^*$ when $V=\mathbf{0}$. This implies $\text{trace}(\tilde{\Sigma}(\mu^*))\le\text{trace}(\Sigma(\mu^*))+\text{trace}(E(\mu^*))$. Since $\mu^*$ is the optimal solution to the priori GKFSP problem when $V=\mathbf{0}$, we have $\text{trace}(\Sigma(\mu^*))\le\text{trace}(\Sigma(\tilde{\mu}_1^*))$. Moreover, one can show that the error covariance of the Kalman filter is always lower bounded (in the positive semi-definite sense) by the error covariance of the Kalman filter with zero measurement noise covariance (with the other system matrices fixed). We obtain $\text{trace}(\Sigma(\tilde{\mu}_1^*))\le\text{trace}(\tilde{\Sigma}(\tilde{\mu}_1^*))$. It then follows from the above arguments that $\text{trace}(\tilde{\Sigma}(\mu^*))\le\text{trace}(\tilde{\Sigma}(\tilde{\mu}_1^*))+\text{trace}(E(\mu^*))$. Similarly, we can show that  $\text{trace}(\tilde{\Sigma}^*(\mu^*))\le\text{trace}(\tilde{\Sigma}^*(\tilde{\mu}_2^*))+\text{trace}(E^*(\mu^*))$.
\end{IEEEproof}
 
The above result has the following interpretation. Consider an instance of the priori (resp., posteriori) GKFSP problem with $V=\tilde{V}\in\mathbb{S}_{+}^n$. If we simply take $V=\mathbf{0}$ and apply the algorithm described in Section \ref{sec:optimal solution to GKFSP}, we will  obtain an optimal solution, denoted as $\mu^*$, to the corresponding instance of the priori (resp., posteriori) GKFSP problem (with $V=\mathbf{0}$). Theorem $\ref{thm:noisy measurement case}$ shows that the performance (i.e., suboptimality) of this sensor placement (i.e., $\mu^*$) for the original priori (resp., posteriori) GKFSP instance with $V=\tilde{V}$ can be bounded by $\text{trace}(\tilde{\Sigma}(\mu^*))\le\text{trace}(\tilde{\Sigma}(\tilde{\mu}_1^*))+\text{trace}(E(\mu^*))$ (resp., $\text{trace}(\tilde{\Sigma}^*(\mu^*))\le\text{trace}(\tilde{\Sigma}^*(\tilde{\mu}_2^*))+\text{trace}(E^*(\mu^*))$), where $\tilde{\mu}_1^*$ (resp., $\tilde{\mu}_2^*$) is the optimal solution to the instance of the priori (resp., posteriori) GKFSP problem when $V=\tilde{V}$. Moreover,  we have from Eq. \eqref{eqn:extra term on priori upper bound} that as $\tilde{V}$ goes to zero, $\text{trace}(E(\mu^*))$ (resp., $\text{trace}(E^*(\mu^*))$) will go to zero, which implies that $\text{trace}(\tilde{\Sigma}(\mu^*))$ (resp., $\text{trace}(\tilde{\Sigma}^*(\mu^*))$) will go to $\text{trace}(\tilde{\Sigma}(\tilde{\mu}_1^*))$ (resp., $\text{trace}(\tilde{\Sigma}^*(\tilde{\mu}_2^*))$). Similar performance bounds can be obtained for the GKFSA and RGKFSP problems, respectively. 

\begin{figure}[htbp]
  \vspace{-0.8cm}
 \center 
 \subfloat[a][GKFSP: True difference vs. Bound]{
 \includegraphics[width=0.5\linewidth]{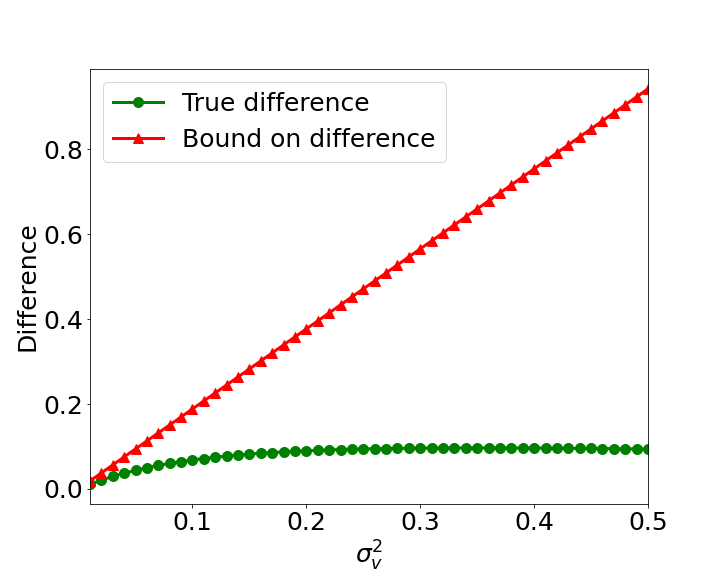}
 }
 \subfloat[b][GKFSP: Suboptimality]{
 \includegraphics[width=0.5\linewidth]{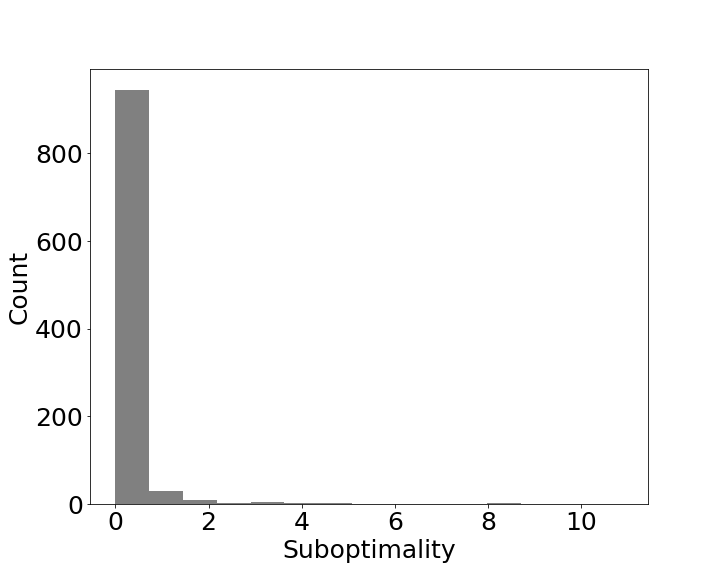}}\\
 \center
 \vspace{-0.55cm}
 \subfloat[c][GKFSA: True difference vs. Bound]{
 \includegraphics[width=0.5\linewidth]{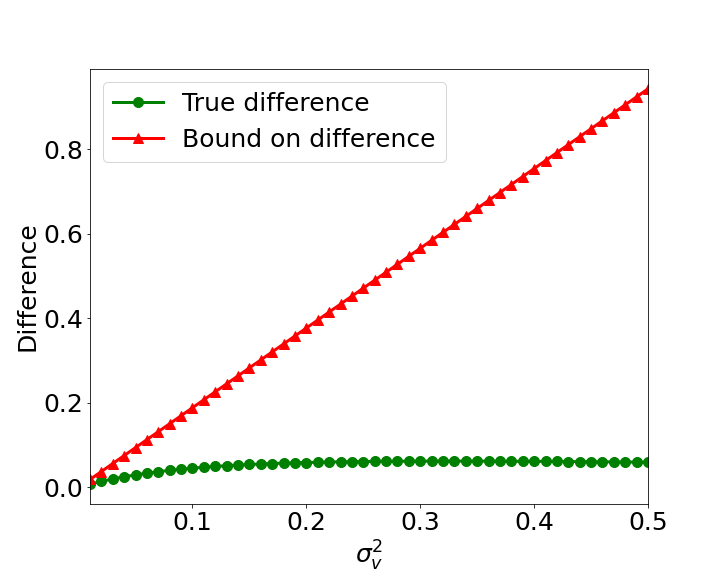}}
 \subfloat[d][GKFSA: Suboptimality]{
 \includegraphics[width=0.5\linewidth]{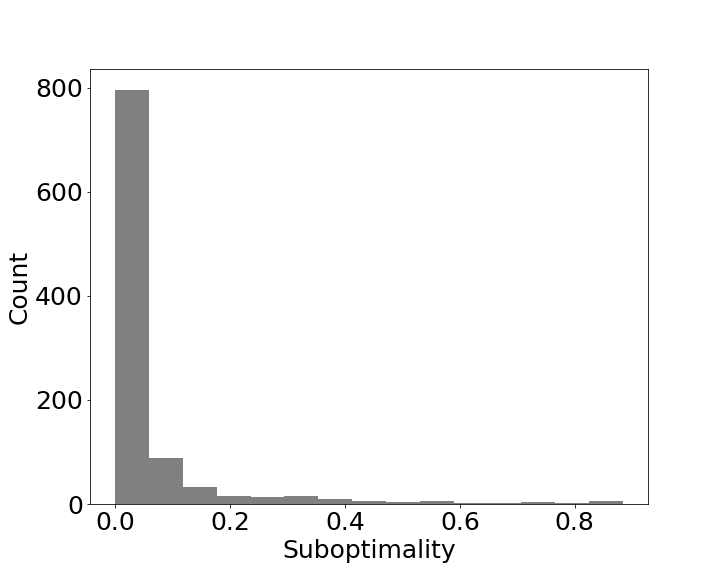}}\\
  \vspace{-0.35cm}
 \subfloat[e][RGKFSP: True difference vs. Bound]{
 \includegraphics[width=0.5\linewidth]{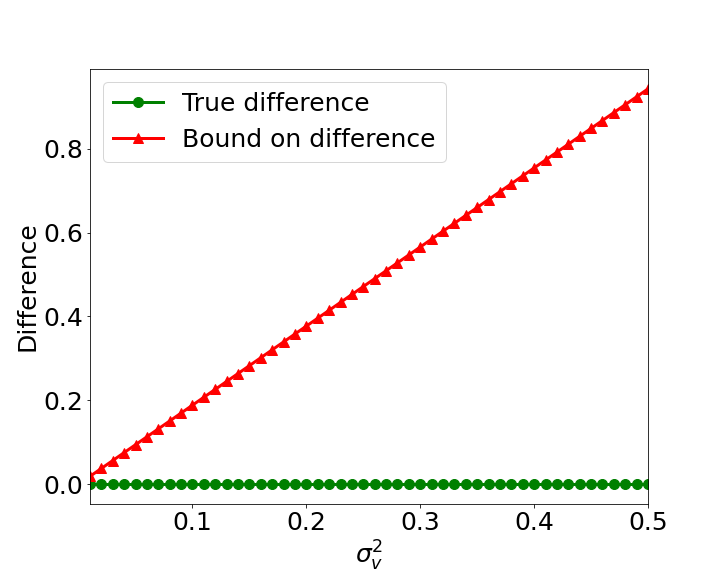}}
 \subfloat[f][RGKFSP: Suboptimality]{
 \includegraphics[width=0.5\linewidth]{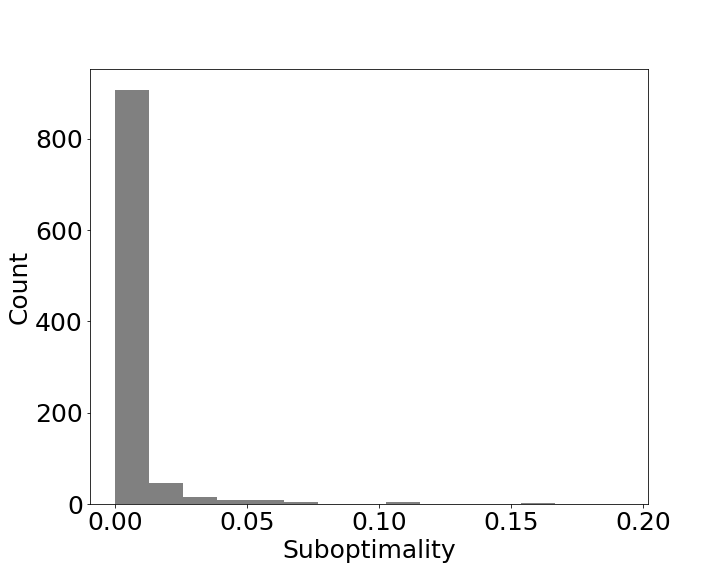}}
 \caption{Performance of the algorithms}
 \label{fig:noise case}
 \end{figure}

We provide simulations to show the performance of the algorithms in Section \ref{sec:optimal solution to GKFSP}, Section \ref{sec:solution to GKFSA}, and Section \ref{sec:alg for RGKFSP}, when applied to solve the GKFSP, GKFSA, and RGKFSP problems with  measurement noise, respectively. Specifically, consider a strongly connected graph $\mathcal{G}(A)$ with $\mathcal{X}(A)=\{x_1,\dots,x_{10}\}$ and $|\mathcal{E}(A)|=15$, where node $x_1$ has the stochastic input with variance $\sigma_w^2=0.1$. Set the measurement matrix $C=I_{10}$ and the sensor noise covariance $V=\sigma^2_v I_{10}$, where $\sigma_v^2\in\mathbb{R}_{\ge0}$. Under a fixed cost $h_i\in\mathbb{Z}_{\ge0}$ to place sensor at $x_i$, a budget $H\in\mathbb{Z}_{\ge0}$, a fixed cost $f_i\in\mathbb{Z}_{\ge0}$ to attack sensor at $x_i$, and an attack budget $F\in\mathbb{Z}_{\ge0}$, we randomly generate the corresponding system dynamics matrix $A\in\mathbb{R}^{10\times 10}$ by selecting each nonzero element of $A$ from a standard normal distribution.   Fig. \ref{fig:noise case}(a) and Fig. \ref{fig:noise case}(b) show  the performance of the algorithm described in Section \ref{sec:optimal solution to GKFSP}, when applied to solve the (priori) GKFSP instances with $V=\sigma_v^2I_{10}$. Specifically, Fig. \ref{fig:noise case}(a) is obtained for a single realization of $A$, which compares the gap (i.e., difference) between the optimal solution to the GKFSP problem (found by brute force and denoted as $OPT$) and the solution returned by the algorithm (denoted as $ALG$), with the bound (on the difference) provided in Theorem $\ref{thm:noisy measurement case}$, when $\sigma_v^2$ ranges from $0.01$ to $0.5$. Fig. \ref{fig:noise case}(b) shows a histogram of the suboptimality of the algorithm, computed as $\frac{ALG-OPT}{OPT}$, over $1000$ realizations of $A$, when $\sigma_v^2=5$. Similarly, Fig. \ref{fig:noise case}(c)-(d) and Fig. \ref{fig:noise case}(e)-(f) show the performance of the algorithm described in Section  \ref{sec:solution to GKFSA} for GKFSA and Algorithm $\ref{algorithm:DP for RGKFSP}$ for RGKFSP, respectively. Note that we fix a sensor placement $\mu$ when solving the GKFSA instances. Moreover, the objective function of RGKFSP associated with the solution returned by Algorithm $\ref{algorithm:DP for RGKFSP}$ is computed against the corresponding optimal sensor attack when $V=\sigma_v^2 I_{10}$. The simulations show that the bounds in Theorem $\ref{thm:noisy measurement case}$ are conservative and that the algorithms (for zero sensor noise) give solutions that are close to optimal for the noisy measurement instances, particularly for RGKSP,  even when $\sigma_w^2/ \sigma_v^2$ becomes small.
 \setkeys{Gin}{width=0.23\textwidth}
\vspace{-0.5cm}

\section{Conclusion}
We considered networked dynamical systems affected by a stochastic input. Under this setting, we first studied the problem for a system designer to optimally place sensors over the network subject to a budget constraint in order to minimize the trace of the steady-state error covariance of the corresponding Kalman filter. We then studied the optimal sensor attack problem where an adversary can attack the placed sensors under an attack budget constraint in order to maximize the trace of the steady-state error covariance of the Kalman filter corresponding to the surviving sensors. Using the graph structure of the networked system, we provided polynomial-time algorithms to solve these two problems. Furthermore, we studied the resilient sensor placement for the system designer when faced with an adversary. We showed that this problem is NP-hard, and provided a pseudo-polynomial-time algorithm to solve it. Although these results are obtained when there is no sensor noise, we provided bounds on the suboptimality of the proposed strategies in the presence of sensor measurement noise. Future work on charactering optimal solutions when there is sensor noise and providing algorithms for systems with multiple stochastic inputs are of interest.
\vspace{-0.1cm}
\section*{Appendix}
\subsection*{Proof of Lemma $\ref{lemma:small noise case}$:}
Denote $\Sigma_{k/k-1}$ (resp., $\Sigma_{k/k}$) as the {\it a priori} (resp., {\it a posteriori}) error covariance of the Kalman filter at time step $k$ when $V=\mathbf{0}$, and denote  $\tilde{\Sigma}_{k/k-1}$ (resp., $\tilde{\Sigma}_{k/k}$) as the {\it a priori} (resp., {\it a posteriori}) error covariance of the Kalman filter at time step $k$ when $V=\tilde{V}$. Denoting $\bar{W}\triangleq BWB^T$, we have (from \cite{anderson1979optimal}):
\begin{equation*}
\tilde{\Sigma}_{k+1/k}=(A-\tilde{K}_kC)\tilde{\Sigma}_{k/k-1}(A-\tilde{K}_kC)^T+\bar{W}+\tilde{K}_k\tilde{V}\tilde{K}_k^T,
\end{equation*}
where $k\ge0$ and $\tilde{K}_k\triangleq A\tilde{\Sigma}_{k/k-1} C^T(C\tilde{\Sigma}_{k/k-1} C^T+\tilde{V})^{-1}$ is the corresponding Kalman gain at time step $k$. For any time step $k$, the Kalman gain $\tilde{K}_k$ satisfies
\begin{equation}
\label{eqn:optimality of kalman gain}
\tilde{K}_k=\mathop{\arg\min}_J\big\{(A-JC)\tilde{\Sigma}_{k/k-1}(A-JC)^T+\bar{W}+J\tilde{V}J^T\big\},
\end{equation}
where the minimization is in the positive semi-definite sense \cite{anderson1979optimal}. Since the pair $(A,BW^{1/2})$ (resp., $(A,C)$) is stabilizable (resp., detectable), we know from a more general version of Lemma $\ref{lemma:Anderson optimal filtering graph}$ for general system matrices in \cite{anderson1979optimal} that the limit $\tilde{\Sigma}=\mathop{\lim}_{k\to\infty}\tilde{\Sigma}_{k+1/k}$ exists, and satisfies 
\begin{equation*}
\label{eqn:alternate dare noise}
\tilde{\Sigma}=(A-\tilde{K}C)\tilde{\Sigma}(A-\tilde{K}C)^T+\bar{W}+\tilde{K}\tilde{V}\tilde{K}^T,
\end{equation*}
where $\tilde{K}\triangleq A\tilde{\Sigma} C^T(C\tilde{\Sigma} C^T+\tilde{V})^{-1}$ is the corresponding (steady-state) Kalman gain. Similarly, we have 
\begin{equation}
\label{eqn:noise free steady state alternate}
\Sigma=(A-KC)\Sigma(A-KC)^T+\bar{W},
\end{equation}
where $K\triangleq A\Sigma C^T(C\Sigma C^T)^{-1}$. Noting the optimality of the Kalman gains from Eq. \eqref{eqn:optimality of kalman gain}, there exists, as argued in \cite{anderson1979optimal}, a suboptimal filter (when $V=\tilde{V}$) with a (time-invariant) suboptimal gain given by $K$ such that the corresponding {\it a priori} error covariance at time step $k+1$, denoted as $\hat{\Sigma}_{k+1/k}$, satisfies 
\begin{equation}
\label{eqn:suboptimal covariance priori}
\hat{\Sigma}_{k+1/k}=(A-KC)\hat{\Sigma}_{k/k-1}(A-KC)^T+\bar{W}+K\tilde{V}K^T.
\end{equation}
Furthermore, the limit $\hat{\Sigma}\triangleq\mathop{\lim}_{k\to\infty}\hat{\Sigma}_{k+1/k}$ exists and satisfies $\hat{\Sigma}\succeq\tilde{\Sigma}$ \cite{anderson1979optimal}. We then obtain from Eq. \eqref{eqn:noise free steady state alternate} and (the steady-state version of) Eq. \eqref{eqn:suboptimal covariance priori} the following:
\begin{equation}
\label{eqn:lyapunov function}
E=(A-KC)E(A-KC)^T+K\tilde{V}K^T,
\end{equation}
where $E=\hat{\Sigma}-\Sigma$. Since the matrix $A-KC$ is stable \cite{anderson1979optimal}, we have that there exists a unique finite positive semi-definite matrix $E$ that satisfies Eq. \eqref{eqn:lyapunov function} and can be written as $E=\sum_{m=0}^{\infty}(A-KC)^mK\tilde{V}K^T((A-KC)^T)^m$ (e.g., \cite{anderson1979optimal}). It then follows from the arguments above that $\hat{\Sigma}=E+\Sigma\succeq\tilde{\Sigma}$.

Similarly, we have from \cite{anderson1979optimal} that $\tilde{\Sigma}_{k/k}$ satisfies
$\tilde{\Sigma}_{k/k}=(I_n-\tilde{L}_kC)\tilde{\Sigma}_{k/k-1}$, where $\tilde{L}_{k}\triangleq\tilde{\Sigma}_{k/k-1}C^T(C\tilde{\Sigma}_{k/k-1}C^T+\tilde{V})^{-1}$. Moreover, the limits $\tilde{\Sigma}^*\triangleq\mathop{\lim}_{k\to\infty}\tilde{\Sigma}_{k/k}$ and $\Sigma^*\triangleq\mathop{\lim}_{k\to\infty}\Sigma_{k/k}$ exist and satisfy $\tilde{\Sigma}^*=(I_n-\tilde{L}C)\tilde{\Sigma}$ and $\Sigma^*=(I_n-LC)\Sigma$, respectively, where $\tilde{L}\triangleq\tilde{\Sigma}C^T(C\tilde{\Sigma}C^T+\tilde{V})^{-1}$ and $L\triangleq\Sigma C^T(C\Sigma C^T)^{-1}$. Similarly, the {\it a posteriori} error covariance at time step $k$ of the suboptimal filter (when $V=\tilde{V}$) as described above, denoted as $\hat{\Sigma}_{k/k}$, is given by
\begin{equation}
\label{eqn:suboptimal covariance posteriori}
\hat{\Sigma}_{k/k}=(I_n-LC)\hat{\Sigma}_{k/k-1}.
\end{equation}
Since the limit $\hat{\Sigma}=\mathop{\lim}_{k\to\infty}\hat{\Sigma}_{k+1/k}$ exists, we know from Eq. \eqref{eqn:suboptimal covariance posteriori} that the limit $\hat{\Sigma}^*\triangleq\mathop{\lim}_{k\to\infty}\hat{\Sigma}_{k/k}$ also exists. Using similar arguments to those in \cite{anderson1979optimal}, one can show that $\hat{\Sigma}^*\succeq\tilde{\Sigma}^*$. Thus, we have $\hat{\Sigma}^*-\Sigma^*=(I_n-LC)(\hat{\Sigma}-\Sigma)=(I_n-LC)E$, which implies $\tilde{\Sigma}^*\preceq\Sigma^*+(I_n-LC)E$.\hfill\IEEEQED
\vspace{-0.2cm}


\bibliographystyle{IEEEtran}
\bibliography{main}
\end{document}